\documentclass[a4paper, final, 12pt]{article}

\usepackage{lipsum}
\usepackage{amsmath, amssymb, latexsym, amscd,amsfonts,amstext}
\usepackage[mathscr]{eucal}
\usepackage{graphicx}
\usepackage{subfig}
\usepackage{float}
\usepackage{color}
\usepackage{hyperref}
\usepackage[utf8]{inputenc}
\usepackage[english]{babel}
\usepackage{pgfplotstable}
\usepackage{tikz}
\usepackage{epstopdf}
\usepackage{algorithmic}
\usepackage{amsopn}

\ifpdf
  \DeclareGraphicsExtensions{.eps,.pdf,.png,.jpg}
\else
  \DeclareGraphicsExtensions{.eps}
\fi

\numberwithin{theorem}{section}
\numberwithin{equation}{section}

\begin{document}

\title{Globally strictly convex cost functional for a 1-D inverse medium scattering problem with experimental data\thanks{ The work of MVK and AEK was supported by US Army Research Laboratory and US Army Research Office grant W911NF-15-1-0233 and by the Office of Naval Research grant N00014-15-1-2330. }}
\author{Michael V. Klibanov \thanks{The corresponding author} \thanks{Department of Mathematics \& Statistics, University of North Carolina at Charlotte, Charlotte, NC 28223, USA},  
Aleksandr E.\ Kolesov \footnotemark[3] \thanks{Institute of Mathematics and Information Science, North-Eastern Federal University, Yakutsk, Russia}, \\ 
Lam Nguyen\thanks{US Army Research Laboratory, 2800 Powder Mill Road Adelphi, MD 20783-1197, USA },  and 
Anders Sullivan\footnotemark[5]}
%\date{}
\maketitle

\begin{abstract}
A new numerical method is proposed for a 1-D inverse medium scattering
problem with multi-frequency data. This method is based on the construction
of a weighted cost functional. The weight is a Carleman Weight Function
(CWF). In other words, this is the function, which is present in the
Carleman estimate for the undelying differential operator. The presence of
the CWF makes this functional strictly convex on any a priori chosen ball
with the center at $\left\{ 0\right\} $ in an appropriate Hilbert space.
Convergence of the gradient minimization method to the exact solution
starting from any point of that ball is proven. Computational results for
both computationally simulated and experimental data show a good accuracy of
this method.
\end{abstract}

%\title{Globally strictly convex cost functional for a 1-D inverse medium scattering problem with experimental data\thanks{ The work of MVK and AEK was supported by US Army Research Laboratory and US Army Research Office grant W911NF-15-1-0233 and by the Office of Naval Research grant N00014-15-1-2330. }}
%\author{Michael V. Klibanov$^{\ast }$\thanks{ $^{\ast }$Department of Mathematics \& Statistics, University of North Carolina at Charlotte, Charlotte, NC 28223, USA}\thanks{The corresponding author}, Aleksandr E.\ Kolesov $^{\ast }$\thanks{$^{\ast }$Department of Mathematics \& Statistics, University of North Carolina at Charlotte, Charlotte, NC 28223, USA }, Lam Nguyen$^{\times }$  \thanks{$^{\times }$US Army Research Laboratory, 2800 Powder Mill Road Adelphy, MD 20783-1197, USA } and Anders Sullivan$^{\times }$ \thanks{ $^{\times }$US Army Research Laboratory, 2800 Powder Mill Road Adelphy, MD 20783-1197, USA }}
%\date{}

%
%\author{Michael V. Klibanov\thanks{
%Department of Mathematics and Statistics, University of North Carolina at
%Charlotte, Charlotte, NC 28223, USA; 
%(\texttt{mklibanv@uncc.edu}, \texttt{hliu34@uncc.edu}, \texttt{dnguye70@uncc.edu}, \texttt{lnguye50@uncc.edu}) }
%\and Hui Liu\footnotemark[2] \and Dinh-Liem Nguyen\footnotemark[2]  \and 
%Loc. H. Nguyen\footnotemark[2]   }

%\emails{mklibanv@uncc.edu, hliu34@uncc.edu, dnguye70@uncc.edu, lnguye50@uncc.edu}



\section{Introduction}

\label{sec:1}

The experimental data used in this paper were collected by the Forward
Looking Radar of the US\ Army Research Laboratory \cite{Radar}. That radar
was built for detection and possible identification of shallow
explosive-like targets. Since targets are three dimensional objects, one
needs to measure a three dimensional information about each target. However,
the radar measures only one time dependent curve for each target, see Figure
5. Therefore, one can hope to reconstruct only a very limited information
about each target. So, we reconstruct only an estimate of the dielectric
constant of each target. For each target, our estimate likely provides a
sort of an average of values of its spatially distributed dielectric
constant. But even this information can be potentially very useful for
engineers. Indeed, currently the radar community is relying only on the
energy information of radar images, see, e.g. \cite{Soumekh}. Estimates of
dielectric constants of targets, if taken alone, cannot improve the current
false alarm rate. However, these estimates can be potentially used as an
additional piece of information. Being combined with the currently used
energy information, this piece of the information might result in the future in new
classification algorithms, which might improve the current false alarm rate.

An Inverse Medium Scattering Problem (IMSP) is often also called a
Coefficient Inverse Problem (CIP). IMSPs/CIPs are both ill-posed and highly
nonlinear. Therefore, an important question to address in a numerical
treatment of such a problem is: \emph{How to reach a sufficiently small
neighborhood of the exact coefficient without any advanced knowledge of this
neighborhood?} The size of this neighborhood should depend only on the level
of noise in the data and on approximation errors. We call a numerical
method, which has a rigorous guarantee of achieving this goal, \emph{%
globally convergent method }(GCM).

In this paper we develop analytically a new globally convergent method for a
1-D Inverse Medium Scattering Problem (IMSP) with the data generated by
multiple frequencies. In addition to the analytical study, we test this
method numerically using both computationally simulated and the above
mentioned experimental data. 

First, we derive a nonlinear integro-differential equation in which the
unknown coefficient is not present. The \emph{new} \emph{element} of this
paper is the method of the solution of this equation. This method is based
on the construction of a weighted least squares cost functional. The key
point of this functional is the presence of the Carleman Weight Function
(CWF) in it. This is the function, which is involved in the Carleman
estimate for the underlying differential operator. We prove that, given a
closed ball of an arbitrary radius $R>0$ with the center at $\left\{
0\right\} $ in an appropriate Hilbert space, one can choose the parameter $%
\lambda >0$ of the CWF in such a way that this functional becomes strictly
convex on that ball.

The existence of the unique minimizer on that closed ball as well as
convergence of minimizers to the exact solution when the level of noise in
the data tends to zero are proven. In addition, \ it is proven that the
gradient projection method reaches a sufficiently small neighborhood of the
exact coefficient if its starting point is an arbitrary point of that ball.
The size of that neighborhood is proportional to the level of noise in the
data. Therefore, since restrictions on $R$ are not imposed in our method,
then this is a \emph{globally convergent} numerical method. We note that in
the conventional case of a non convex cost functional a gradient-like method
converges to the exact solution only if its starting point is located in a
sufficiently small neighborhood of this solution: this is due to the
phenomenon of multiple local minima and ravines of such functionals.

Unlike previously developed globally convergent numerical methods of the
first type for CIPs (see this section below), the convergence analysis for
the technique of the current paper does not impose a smallness condition on
the interval $\left( \underline{k},\overline{k}\right) $ of the variations
of the wave numbers $k\in \left( \underline{k},\overline{k}\right) \subset
\left\{ k>0\right\} $. 

The majority of currently known numerical methods of solutions of nonlinear
ill-posed problems use the nonlinear optimization. In other words, a least
squares cost functional is minimized in each problem, see, e.g. \cite%
{Chavent,Engl,Gonch1,Gonch2}. However, the major problem with these
functionals is that they are usually non convex. Figure 1 of the paper \cite%
{Scales} presents a numerical example of multiple local minima and ravines
of non-convex least squares cost functionals for some CIPs. Hence,
convergence of the optimization process of such a functional to the exact
solution can be guaranteed only if a good approximation for that solution is
known in advance. However, such an approximation is rarely available in
applications. This prompts the development of globally convergent numerical
methods for CIPs, see, e.g. \cite{BKSISC, BK, BK2, Chow, Karch, Klib95, Klib97, KT,  KK, Klibfreq,  KlibLoc, KTSIAP, Kuzh, IEEE,  Exp2, TBKF2}%
.

The first author with coauthors has proposed two types of GCM for CIPs with
single measurement data. The GCM of the first type is reasonable to call the
\textquotedblleft tail functions method". This development has started from
the work \cite{BKSISC} and has been continued since then, see, e.g. \cite{BK, Chow,  Klibfreq, KlibLoc, Kuzh, IEEE, Exp2, TBKF2} and references cited
therein. In this case, on each step of an iterative process one solves the
Dirichlet boundary value problem for a certain linear elliptic PDE, which
depends on that iterative step. The solution of this PDE allows one to
update the unknown coefficient first and then to update a certain function,
which is called \textquotedblleft the tail function". The convergence
theorems for this method impose a smallness condition on the interval of the
variation of either the parameter $s>0$ of the Laplace transform of the
solution of a hyperbolic equation or of the wave number $k>0$ in the
Helmholtz equation. Recall that the method of this paper does not impose the
latter assumption.

In this paper we present a new version of the GCM of the second type. In any
version of the GCM of the second type a weighted cost functional with a CWF
in it is constructed. The same properties of the global strict convexity and
the global convergence of the gradient projection method hold as the ones
indicated above. The GCM of the second type was initiated in \cite%
{Klib95,Klib97,KT} with a recently renewed interest in \cite{BK2,KK,KTSIAP}.
The idea of any version of the GCM of the second type has direct roots in
the method of \cite{BukhKlib}, which is based on Carleman estimates and
which was originally designed in \cite{BukhKlib} only for proofs of
uniqueness theorems for CIPs, also see the recent survey in \cite{Ksurvey}.

Another version of the GCM with a CWF in it was recently developed in \cite%
{Bau1} for a CIP for the hyperbolic equation $w_{tt}=\Delta w+a\left(
x\right) w+f\left( x,t\right) ,$ where $a\left( x\right) $ is the unknown
coefficient. This GCM was tested numerically in \cite{Bau2}. In \cite%
{Bau1,Bau2} non-vanishing conditions are imposed: it is assumed that either $%
f\left( x,0\right) \neq 0$ or $w\left( x,0\right) \neq 0$ or $w_{t}\left(
x,0\right) \neq 0$ in the entire domain of interest. Similar assumptions are
imposed in \cite{BK2,KK} for the GCM of the second type. On the other hand,
we consider in the current paper, so as in \cite{Klib95,Klib97,KT,KTSIAP},
the fundamental solution of the corresponding PDE. The differences between
the fundamental solutions of those PDEs and solutions satisfying
non-vanishing conditions cause quite significant differences between \cite%
{Klib95,Klib97,KT,KTSIAP} and \cite{Bau1,Bau2,BK2,KK} of corresponding
versions of the GCM of the second type.

Recently, the idea of the GCM of the second type was extended to the case of
ill-posed Cauchy problems for quasilinear PDEs, see the theory in \cite%
{Klquasi} and some extensions and numerical examples in \cite%
{BakKlKosh,KlKosh}.

CIPs of wave propagation are a part of a bigger subfield, Inverse Scattering
Problems (ISPs). ISPs attract a significant attention of the scientific
community. In this regard we refer to some direct methods which successfully
reconstruct positions, sizes and shapes of scatterers without iterations 
\cite{CakoniColton:gm2003, CharalambopoulosGintidesKiriaki:ip2002,
Ito,KirschGrinberg:ols2008, LiLiuZou:smms2014, LiuSini:sjsc2009,
Liem:ipi2016,Sini}. We also refer to \cite{AmmariKang:lnim2004,
LiuSini:sjsc2009, Novikov1992, Novikov:faa1988} for some other ISPs in the
frequency domain. In addition, we cite some other numerical methods for ISPs
considered in \cite{AmmariChowZou,Bao:ip2015, Sini2}.

As to the CIPs with multiple measurement, i.e. the Dirichlet-to-Neumann map
data, we mention recent works \cite{Agal,Kab,Novikov} and references cited
therein, where reconstruction procedures are developed, which do not require
a priori knowledge of a small neighborhood of the exact coefficient.

In section 2 we state our inverse problem. In section 3 we construct that
weighted cost functional. In section 4 we prove the main property of this
functional: its global strict convexity. In section 5 we prove the global
convergence of the gradient projection method of the minimization of this
functional. Although this paper is mostly an analytical one (sections 3-5),
we complement the theory with computations. In section 6 we test our method
on computationally simulated data. In section 7 we test it on experimental
data. Concluding remarks are in section 8.

\section{Problem statement}

\label{sec:2}

\subsection{Statement of the inverse problem}

\label{sec:2.1}

Let the function $c\left( x\right) ,x\in \mathbb{R}$ be the spatially
distributed dielectric constant of the medium. We assume that%
\begin{equation}
c\in C^{2}\left( \mathbb{R}\right) ,c\left( x\right) \geq 1,\forall x\in 
\mathbb{R},  \label{2.1}
\end{equation}%
\begin{equation}
c\left( x\right) =1,\forall x\notin \left( 0,1\right) .  \label{2.2}
\end{equation}%
Fix the source position $x_{0}<0.$ For brevity, we do not indicate below
dependence of our functions on $x_{0}.$ Consider the 1-D Helmholtz equation
for the function $u\left( x,k\right) $,%
\begin{equation}
u^{\prime \prime }+k^{2}c\left( x\right) u=-\delta \left( x-x_{0}\right)
,x\in \mathbb{R},  \label{2.4}
\end{equation}%
\begin{equation}
\lim_{x\rightarrow \infty }\left( u^{\prime }+iku\right)
=0,\lim_{x\rightarrow -\infty }\left( u^{\prime }-iku\right) =0.  \label{2.6}
\end{equation}%
Let $u_{0}\left( x,k\right) $ be the solution of the problem (\ref{2.4}), (%
\ref{2.6}) for the case $c\left( x\right) \equiv 1.$ Then%
\begin{equation}
u_{0}\left( x,k\right) =\frac{\exp \left( -ik\left\vert x-x_{0}\right\vert
\right) }{2ik}.  \label{2.60}
\end{equation}%
Our interest is in the following inverse problem:

\textbf{Inverse Medium Scattering Problem (IMSP)}.\emph{\ Let }$[\underline{k%
},\overline{k}]\subset \left( 0,\infty \right) $\emph{\ be an interval of
wavenumbers }$k$\emph{. Reconstruct the function }$c\left( x\right) ,$\emph{%
\ assuming that the following function }$g_{0}\left( k\right) $\emph{\ is
known} 
\begin{equation}
g_{0}\left( k\right) =\frac{u(0,k)}{u_{0}(0,k)},k\in \lbrack \underline{k},%
\overline{k}].  \label{2.8}
\end{equation}

Denote%
\begin{equation}
w\left( x,k\right) =\frac{u\left( x,k\right) }{u_{0}\left( x,k\right) }.
\label{2.100}
\end{equation}%
It follows from (\ref{2.8}), (\ref{2.100}) and \cite{KlibLoc} that 
\begin{equation}
w\left( 0,k\right) =g_{0}\left( k\right) ,k\in \lbrack \underline{k},%
\overline{k}],  \label{2.101}
\end{equation}%
\begin{equation}
w^{\prime }\left( 0,k\right) =g_{1}\left( k\right) =2ik\left( g_{0}\left(
k\right) -1\right) ,k\in \lbrack \underline{k},\overline{k}].  \label{2.160}
\end{equation}

\subsection{Some properties of the solution of forward and inverse problems}

\label{sec:2.2}

In this subsection we briefly outline some results of \cite{KlibLoc}, which
we use below in this paper. Existence and uniqueness of the solution $%
u\left( x,k\right) \in C^{3}\left( \mathbb{R}\right) $ for each $k>0$ was
established in \cite{KlibLoc}. Also, it was proven in \cite{KlibLoc} that 
\begin{equation}
u\left( x,k\right) \neq 0,\forall x\in \left[ 0,1\right] ,\forall k>0.
\label{2.9}
\end{equation}%
In particular, $g_{0}\left( k\right) \neq 0,\forall k\in \lbrack \underline{k%
},\overline{k}].$ In addition, uniqueness of our IMSP was proven in \cite%
{KlibLoc}. Also, the following asymptotic behavior of the function $u\left(
x,k\right) $ takes place: 
\begin{equation}
u\left( x,k\right) =\frac{1}{2ikc^{1/4}\left( x\right) }\exp \left[
-ik\int\limits_{x_{0}}^{x}\sqrt{c\left( \xi \right) }d\xi \right] \left( 1+%
\widehat{u}\left( x,k\right) \right) ,k\rightarrow \infty ,\forall x\in %
\left[ 0,1\right] ,  \label{2.10}
\end{equation}%
\begin{equation}
\widehat{u}\left( x,k\right) =O\left( \frac{1}{k}\right) ,\partial _{k}%
\widehat{u}\left( x,k\right) =O\left( \frac{1}{k^{2}}\right) ,k\rightarrow
\infty .  \label{2.1000}
\end{equation}

Given (\ref{2.9}) and (\ref{2.10}) we now can uniquely define the function $%
\log w\left( x,k\right) $ as in \cite{KlibLoc}. The difficulty here is in
defining $\mathop{\rm Im}\left( \log u\left( x,k\right) \right) ,$ since
this number is usually defined up to the addition of $2n\pi ,$ where $n$ is
an integer. For sufficiently large values of $k$ we define the function $%
\log w\left( x,k\right) $ using (\ref{2.60}), (\ref{2.100}), (\ref{2.10})
and (\ref{2.1000}) as 
\begin{equation}
\log w\left( x,k\right) =-\frac{1}{4}\ln c\left( x\right) -ik\left(
\int\limits_{x_{0}}^{x}\sqrt{c\left( \xi \right) }d\xi -x+x_{0}\right) +%
\widehat{w}\left( x,k\right) ,x\in \left( 0,1\right) ,  \label{2.11}
\end{equation}%
where 
\begin{equation}
\widehat{w}\left( x,k\right) =O\left( \frac{1}{k}\right) ,\partial _{k}%
\widehat{w}\left( x,k\right) =O\left( \frac{1}{k^{2}}\right) ,k\rightarrow
\infty .  \label{2.110}
\end{equation}%
Hence, for sufficiently large $k$, 
\begin{equation}
\left\vert \widehat{w}\left( x,k\right) \right\vert <2\pi ,  \label{2.12}
\end{equation}%
which eliminates the above mentioned ambiguity. Suppose that the number $%
\overline{k}$ is so large that (\ref{2.12}) is true for $k\geq \overline{k}.$
Then $\log w\left( x,\overline{k}\right) $ is defined as in (\ref{2.11}). As
to not large values of $k$, we define the function (\ref{2.11})$=\log w(x,k)$
as 
\begin{equation}
\psi (x,k)=-\int\limits_{k}^{\overline{k}}\frac{\partial _{k}w(x,\xi )}{%
w(x,\xi )}d\xi +\log w(x,\overline{k}).  \label{2.13}
\end{equation}%
By (\ref{2.9}) $w(x,\xi )\neq 0,\forall x\in \left[ 0,1\right] ,\forall \xi
>0.$ Differentiating both sides of (\ref{2.13}) with respect to $k$, we
obtain 
\begin{equation}
\partial _{k}w(x,k)-w(x,k)\partial _{k}\psi (x,k)=0.  \label{2.14}
\end{equation}%
Multiplying both sides of (\ref{2.14}) by $\exp (-\psi (x,k))$, we obtain $%
\partial _{k}\left( e^{-\psi (x,k)}w(x,k)\right) =0.$ Hence, there exists a
function $C=C\left( x\right) $ independent on $k$ such that 
\begin{equation}
w(x,k)=C\left( x\right) e^{\psi (x,k)}.  \label{2.15}
\end{equation}%
Setting in (\ref{2.15}) $k=\overline{k}$ and using the fact that by (\ref%
{2.13}) $\psi (x,\overline{k})=\log w(x,\overline{k})$, we obtain 
\begin{equation}
C=C\left( x\right) =1,x\in \left[ 0,1\right] .  \label{2.150}
\end{equation}%
Hence, (\ref{2.13}) and (\ref{2.15}) imply that $\log w(x,k)$ is defined as $%
\log w(x,k)=\psi (x,k).$

\section{The Weighted Cost Functional}

\label{sec:3}

In this section we construct the above mentioned weighted cost functional
with the CWF in it.

\textbf{Lemma 3.1} (Carleman estimate). \emph{For any complex valued
function }$u\in H^{2}\left( 0,1\right) $\emph{\ with }$u(0)=u^{\prime }(0)=0$%
\emph{\ and for any parameter }$\lambda >1$\emph{\ the following Carleman
estimate holds }%
\begin{equation}
\int\limits_{0}^{1}\left\vert u^{\prime \prime }\right\vert ^{2}e^{-2\lambda
x}dx\geq C\left[ \int\limits_{0}^{1}|u^{\prime \prime }|^{2}e^{-2\lambda
x}dx+\lambda \int\limits_{0}^{1}|u^{\prime }|^{2}e^{-2\lambda x}dx+\lambda
^{3}\int\limits_{0}^{1}|u|^{2}e^{-2\lambda x}dx\right] ,  \label{3.00}
\end{equation}%
\emph{where\ the constant }$C>0$ \emph{is} \emph{independent of }$u$\emph{\
and }$\lambda .$

\textbf{Proof}. In the case when the integral with $u^{\prime \prime }$ is
absent in the right hand side of (\ref{3.00}) this lemma was proved in \cite%
{KlibLoc}. To incorporate this integral, we note that 
\begin{equation}
2\int\limits_{0}^{1}\left\vert u^{\prime \prime }\right\vert
^{2}e^{-2\lambda x}dx\geq \int\limits_{0}^{1}\left\vert u^{\prime \prime
}\right\vert ^{2}e^{-2\lambda x}dx+C\left[ \lambda
\int\limits_{0}^{1}|u^{\prime }|^{2}e^{-2\lambda x}dx+\lambda
^{3}\int\limits_{0}^{1}|u|^{2}e^{-2\lambda x}dx\right] .  \label{3.02}
\end{equation}%
Let $\widetilde{C}=\min \left( C/2,1/2\right) .$ Then (\ref{3.02}) implies (%
\ref{3.00}) where $C$ is replaced with $\widetilde{C}.$ $\square $

\subsection{Nonlinear integro-differential equation}

\label{sec:3.1}

For $x\in \lbrack 0,1],k\in \lbrack \underline{k},\overline{k}]$ consider
the function $v(x,k)$ and its $k-$derivative $q\left( x,k\right) $, where 
\begin{equation}
v(x,k)=\frac{\log w(x,k)}{k^{2}},\text{ }q\left( x,k\right) =\partial
_{k}v\left( x,k\right) .  \label{3.0}
\end{equation}

Hence,%
\begin{equation}
v\left( x,k\right) =-\int\limits_{k}^{\overline{k}}q\left( x,\tau \right)
d\tau +v\left( x,\overline{k}\right) .  \label{3.1}
\end{equation}%
Consider the function $V\left( x\right) =V\left( x,\overline{k}\right) $,
which we call the \textquotedblleft tail function", and this function is
unknown,%
\begin{equation}
V\left( x\right) =v\left( x,\overline{k}\right) .  \label{3.2}
\end{equation}

Let $\beta \left( x\right) =c\left( x\right) -1.$ Note that since $c\left(
x\right) =1$ for $x\geq 1,$ then equation (\ref{2.4}) and the first
condition (\ref{2.6}) imply that $u\left( x,k\right) =A\left( k\right)
e^{-ikx}$ for $x\geq 1.$ Hence, (\ref{2.60}) and (\ref{2.100}) imply that $%
w^{\prime }\left( x,k\right) =0$ for $x\geq 1.$ It follows from (\ref{2.4}),
(\ref{2.60}), (\ref{2.100})--(\ref{2.160}), (\ref{2.15}) and (\ref{2.150})
that 
\begin{equation}
w^{\prime \prime }-2ikw^{\prime }+k^{2}\beta \left( x\right) w=0,x\in \left(
0,1\right) ,  \label{3.3}
\end{equation}%
\begin{equation}
w\left( 0,k\right) =g_{0}\left( k\right) ,w^{\prime }\left( 0,k\right)
=g_{1}\left( k\right) ,w^{\prime }\left( 1,k\right) =0.  \label{3.4}
\end{equation}%
Using (\ref{2.15}), (\ref{2.150}), (\ref{3.0}) and (\ref{3.3}), we obtain 
\begin{equation}
v^{\prime \prime }+k^{2}\left( v^{\prime }\right) ^{2}-2ikv^{\prime }=-\beta
\left( x\right) .  \label{3.5}
\end{equation}%
Differentiate (\ref{3.5}) with respect to $k$ and use (\ref{3.0})-(\ref{3.4}%
). We obtain 
\begin{equation*}
q^{\prime \prime }-2ikq^{\prime }+2k^{2}q^{\prime }\left( -\int\limits_{k}^{%
\overline{k}}q^{\prime }\left( x,\tau \right) d\tau +V^{\prime }\left(
x\right) \right) -2i\left( -\int\limits_{k}^{\overline{k}}q^{\prime }\left(
x,\tau \right) d\tau +V^{\prime }\left( x\right) \right)
\end{equation*}%
\begin{equation}
+2k\left( -\int\limits_{k}^{\overline{k}}q^{\prime }\left( x,\tau \right)
d\tau +V^{\prime }\left( x\right) \right) ^{2}=0,\text{ }x\in \left(
0,1\right) ,k\in \left( \underline{k},\overline{k}\right) ,  \label{3.6}
\end{equation}%
\begin{equation}
q\left( 0,k\right) =p_{0}\left( k\right) ,q^{\prime }\left( 0,k\right)
=p_{1}\left( k\right) ,q^{\prime }\left( 1,k\right) =0,k\in \left( 
\underline{k},\overline{k}\right) ,  \label{3.7}
\end{equation}%
where%
\begin{equation}
p_{0}\left( k\right) =\frac{\partial }{\partial k}\left( \frac{\log g_{0}(k)%
}{k^{2}}\right) ,p_{1}\left( k\right) =\frac{\partial }{\partial k}\left[ 
\frac{2i}{k}\left( 1-\frac{1}{g_{0}\left( k\right) }\right) \right] ,k\in %
\left[ \underline{k},\overline{k}\right] .  \label{3.70}
\end{equation}

We have obtained an integro-differential equation (\ref{3.6}) for the
function $q$ with the overdetermined boundary conditions (\ref{3.7}). The
tail function $V\left( x\right) $ is also unknown. First, we will
approximate the tail function $V$. Next, we will solve the problem (\ref{3.6}%
), (\ref{3.7}) for the function $q$. To solve this problem, we will
construct the above mentioned weighted cost functional with the CWF $%
e^{-2\lambda x}$ in it, see (\ref{3.00}). This construction, combined with
corresponding analytical results, is the \emph{central }part of our paper.
Thus, even though the problem (\ref{3.6})-(\ref{3.70}) is the same as the
problem (65), (66) in \cite{KlibLoc}, the numerical method of the solution
of the problem (\ref{3.6})-(\ref{3.70}) is \emph{radically} different from
the one in \cite{KlibLoc}.

Now, suppose that we have obtained approximations for both functions $%
V\left( x\right) $ and $q\left( x,k\right) $. Then we obtain the unknown
coefficient $\beta \left( x\right) $ via backwards calculations. First, we
calculate the approximation for the function $v\left( x,k\right) $ via (\ref%
{3.1}) and (\ref{3.2}).\ Next, we calculate the function $\beta \left(
x\right) =c\left( x\right) -1$ via (\ref{3.5}). We have learned from our
numerical experience that the best value of $k$ to use in (\ref{3.5}) for
the latter calculation is $k=\underline{k}.$

\subsection{Approximation for the tail function $V\left( x\right) $}

\label{sec:3.2}

The approximation for the tail function is done here the same way as the
approximation for the so-called \textquotedblleft first tail function" in
section 4.2 of \cite{KlibLoc}. However, while tail functions are updated in 
\cite{KlibLoc}, we are not doing such updates here.

It follows from (\ref{2.100})-(\ref{2.110}) and (\ref{3.0})-(\ref{3.2}) that
there exists a function $r\left( x\right) \in C^{2}\left[ 0,1\right] $ such
that 
\begin{equation}
V\left( x,k\right) =\frac{r\left( x\right) }{k}+O\left( \frac{1}{k^{2}}%
\right) ,\text{ }q\left( x,k\right) =-\frac{r\left( x\right) }{k^{2}}%
+O\left( \frac{1}{k^{3}}\right) ,k\rightarrow \infty ,x\in \left( 0,1\right)
.  \label{3.8}
\end{equation}%
Hence, assuming that the number $\overline{k}$ is sufficiently large, we
drop terms $O\left( 1/\overline{k}^{2}\right) $ and $O\left( 1/\overline{k}%
^{3}\right) $ in (\ref{3.8}).\ Next, we set%
\begin{equation}
V\left( x,k\right) =\frac{r\left( x\right) }{k},\text{ }q\left( x,k\right) =-%
\frac{r\left( x\right) }{k^{2}},k\geq \overline{k},x\in \left( 0,1\right) .
\label{3.9}
\end{equation}%
Set $k:=\overline{k}$ in (\ref{3.6}) and (\ref{3.7}). Next, substitute (\ref%
{3.9}) in (\ref{3.6}) and (\ref{3.7}) at $k=\overline{k}$. We obtain $%
r^{\prime \prime }=0.$ Recall that functions $g_{0}$ and $g_{1}$ are linked
via (\ref{2.160}). Thus, 
\begin{align}
V^{\prime \prime }& =0,\quad \mbox{in }(0,1),  \label{3.10} \\
V(0)& =\frac{\log g_{0}(\overline{k})}{\overline{k}^{2}},V^{\prime }(0)=%
\frac{i}{\overline{k}}\left( 1-\frac{1}{g_{0}(\overline{k})}\right)
,V^{\prime }(1)=0,  \label{3.11}
\end{align}%
where functions $g_{0}$ and $g_{1}(k)$ are defined in (\ref{2.101}) and (\ref%
{2.160}) respectively. It seems to be at the first glance that one can find
the function \ $V$ as, for example Cauchy problem for ODE (\ref{3.10}) with
data $V(0)$ and $V^{\prime }(0).$ However, it was noticed in Remark 5.1 of 
\cite{Kuzh} that this approach, being applied to a similar problem, does not
lead to good results. We have the same observation in our numerical studies.
This is likely to the approximate nature of (\ref{3.9}). Thus, just like in 
\cite{KlibLoc}, we solve the problem (\ref{3.10}), (\ref{3.11}) by the
Quasi-Reversibility Method (QRM). The boundary condition $V^{\prime }(1)=0$
provides a better stability property. \ 

So, we minimize the following functional $J_{\alpha }(V)$ on the set $W$,
where 
\begin{equation}
J_{\alpha }(V)=\frac{1}{2}\left( \Vert V^{\prime \prime }\Vert
_{L^{2}(0,1)}^{2}+\alpha \Vert V\Vert _{H^{3}(0,1)}^{2}\right) ,
\label{3.12}
\end{equation}%
\begin{equation}
V\in W:=\{V\in H^{3}(0,1):V(0)=\frac{\log g_{0}(\overline{k})}{\overline{k}%
^{2}},V^{\prime }\left( 0\right) =\frac{i}{\overline{k}}\left( 1-\frac{1}{%
g_{0}(\overline{k})}\right) ,V^{\prime }(1)=0\},  \label{3.13}
\end{equation}%
where $\alpha >0$ is the regularization parameter. The existence and
uniqueness of the solution of this minimization problem as well as
convergence of minimizers $V_{\alpha }$ in the $H^{2}\left( 0,1\right) -$%
norm to the exact solution $V^{\ast }$ of the problem (\ref{3.11}), (\ref%
{3.12}) with the exact data $g_{0}^{\ast }(\overline{k})$ as $\alpha
\rightarrow 0$ were proved in \cite{KlibLoc}. We note that in the
regularization theory one always assumes existence of an ideal exact
solution with noiseless data \cite{BK,Engl}.

Recall that by the embedding theorem $H^{2}\left( 0,1\right) \subset C^{1}%
\left[ 0,1\right] $ and 
\begin{equation}
\left\Vert f\right\Vert _{C^{1}\left[ 0,1\right] }\leq C\left\Vert
f\right\Vert _{H^{2}\left( 0,1\right) },\forall f\in H^{2}\left( 0,1\right) ,
\label{3.130}
\end{equation}
where $C>0$ is a generic constant$.$ Theorem 3.1 is a reformulation of
Theorem 4.2 of \cite{KlibLoc}.

\textbf{Theorem \ 3.1.} \emph{Let the function }$c^{\ast }\left( x\right) $%
\emph{\ satisfying conditions (\ref{2.1})-(\ref{2.2}) be the exact solution
of our IMSP with the noiseless data }$g_{0}^{\ast }\left( k\right) =w^{\ast
}\left( 0,k\right) ,k\in \left[ \underline{k},\overline{k}\right] $\emph{,
where }$w^{\ast }\left( x,k\right) =u^{\ast }\left( x,k\right) /u_{0}\left(
x,k\right) $\emph{\ and }$u^{\ast }\left( x,k\right) $\emph{\ is the
solution of the forward problem (\ref{2.4}), (\ref{2.6}). Let the exact tail
function }$V^{\ast }\left( x,\overline{k}\right) =\overline{k}^{-2}\log
w^{\ast }\left( x,\overline{k}\right) $\emph{\ and the function }$q^{\ast
}\left( x,\overline{k}\right) =\partial _{k}V^{\ast }\left( x,k\right) \mid
_{k=\overline{k}}$\emph{have the form (\ref{3.9}) with }$r:=r^{\ast }\left(
x\right) .$\emph{\ Assume that for }$k\in \left[ \underline{k},\overline{k}%
\right] $\emph{\ }%
\begin{equation}
\left\vert \log g_{0}\left( k\right) -\log g_{0}^{\ast }\left( k\right)
\right\vert \leq \delta ,\left\vert g_{0}\left( k\right) -g_{0}^{\ast
}\left( k\right) \right\vert \leq \delta ,\left\vert g_{0}^{\prime }\left(
k\right) -\left( g_{0}^{\ast }\right) ^{\prime }\left( k\right) \right\vert
\leq \delta ,  \label{3.14}
\end{equation}%
\emph{where }$\delta >0$\emph{\ is a sufficiently small number, which
characterizes the level of the error in the boundary data. Let in (\ref{3.12}%
) }$\alpha =\alpha \left( \delta \right) =\delta ^{2}.$\emph{\ Let the
function }$V_{\alpha \left( \delta \right) }\left( x\right) \in H^{3}\left(
0,1\right) $\emph{\ be the minimizer of the functional (\ref{3.12}) on the
set of functions }$W$\emph{\ defined in (\ref{3.13}). Then there exists a
constant }$C_{1}=C_{1}\left( \overline{k},c^{\ast }\right) >0$\emph{\
depending only on }$\overline{k}$\emph{\ and }$c^{\ast }$\emph{\ such that} 
\begin{equation}
\left\Vert V_{\alpha \left( \delta \right) }\left( x\right) -V^{\ast }\left(
x,\overline{k}\right) \right\Vert _{C^{1}\left[ 0,1\right] }\leq C\left\Vert
V_{\alpha \left( \delta \right) }\left( x\right) -V^{\ast }\left( x,%
\overline{k}\right) \right\Vert _{H^{2}\left( 0,1\right) }\leq C_{1}\delta .
\label{3.15}
\end{equation}

\textbf{Remark 3.1}. We have also tried to consider two terms in the
asymptotic expansion for $V$ in (\ref{3.8}): the second one with $1/k^{2}.$
This resulted in a nonlinear system of two equations. We have solved it by
via minimizing an analog of the functional of section 3.3. However, the
quality of resulting images deteriorated as compared with the above function 
$V_{\alpha \left( \delta \right) }\left( x\right) .$ In addition, we have
tried to iterate with respect to the tail function $V$. However, the quality
of resulting images has also deteriorated.

\subsection{The weighted cost functional}

\label{sec:3.3}

Consider the function $q\left( x,k\right) $ satisfying (\ref{3.6})-(\ref%
{3.70}). In sections 5.2 and 5.3 we use Lemma 2.1 and Theorem 2.1 of \cite%
{BakKlKosh}. To apply theorems, we need to have zero boundary conditions at $%
x=0,1.$ Hence, we introduce the function $p\left( x,k\right) ,$%
\begin{equation}
p\left( x,k\right) =q\left( x,k\right) -\left( x^{2}-1\right)
^{2}p_{0}-x\left( x^{2}-1\right) ^{2}p_{1},\text{ where }p_{0}=p_{0}\left(
k\right) ,p_{1}=p_{1}\left( k\right) .  \label{3.16}
\end{equation}%
Denote 
\begin{equation}
m\left( x,k\right) =\left( x^{2}-1\right) ^{2}p_{0}\left( k\right) +x\left(
x^{2}-1\right) ^{2}p_{1}\left( k\right) .  \label{3.170}
\end{equation}%
Also, replace in (\ref{3.6}) $V$ with $V_{\alpha \left( \delta \right) }.$
Then (\ref{3.6}), (\ref{3.7}) and (\ref{3.16}) and (\ref{3.170}) imply that%
\begin{equation*}
L\left( p\right) =p^{\prime \prime }+m^{\prime \prime }-2ik\left( p^{\prime
}+m^{\prime }\right) +2k^{2}\left( p^{\prime }+m^{\prime }\right) \left(
-\int\limits_{k}^{\overline{k}}\left( p^{\prime }+m^{\prime }\right) \left(
x,\tau \right) d\tau +V_{\alpha \left( \delta \right) }^{\prime }\left(
x\right) \right)
\end{equation*}%
\begin{equation}
-2i\left( -\int\limits_{k}^{\overline{k}}\left( p^{\prime }+m^{\prime
}\right) \left( x,\tau \right) d\tau +V^{\prime }\left( x\right) \right)
+2k\left( -\int\limits_{k}^{\overline{k}}\left( p^{\prime }+m^{\prime
}\right) \left( x,\tau \right) d\tau +V_{\alpha \left( \delta \right)
}^{\prime }\left( x\right) \right) ^{2}=0  \label{3.17}
\end{equation}%
\begin{equation}
p\left( 0,k\right) =0,p^{\prime }\left( 0,k\right) =0,p^{\prime }\left(
1,k\right) =0.  \label{3.18}
\end{equation}

Introduce the Hilbert space $H$ of pairs of real valued functions $f\left(
x,k\right) =\left( f_{1}\left( x,k\right) ,f_{2}\left( x,k\right) \right) ,$ 
$\left( x,k\right) \in \left( 0,1\right) \times \left( \underline{k},%
\overline{k}\right) $ as%
\begin{equation}
H=\left\{ 
\begin{array}{c}
f\left( x,k\right) :f\left( 0,k\right) =f^{\prime }\left( 0,k\right)
=f^{\prime }\left( 1,k\right) =0, \\ 
\left\Vert f\right\Vert _{H}=\left[ \displaystyle\int\limits_{\underline{k}%
}^{\overline{k}}\left\Vert f\left( x,k\right) \right\Vert _{H^{2}\left(
0,1\right) }^{2}dk\right] ^{1/2}<\infty%
\end{array}%
\right\} .  \label{3.19}
\end{equation}%
Here and below $\left\Vert f\left( x,k\right) \right\Vert _{H^{2}\left(
0,1\right) }^{2}=\left\Vert f_{1}\left( x,k\right) \right\Vert _{H^{2}\left(
0,1\right) }^{2}+\left\Vert f_{2}\left( x,k\right) \right\Vert _{H^{2}\left(
0,1\right) }^{2}.$

Based on (\ref{3.17}) and (\ref{3.18}), we define our weighted cost
functional as%
\begin{equation}
J_{\lambda }\left( p\right) =e^{2\lambda }\displaystyle\int\limits_{%
\underline{k}}^{\overline{k}}\displaystyle\int\limits_{0}^{1}\left\vert
L\left( p\right) \right\vert ^{2}e^{-2\lambda x}dxdk,\forall p\in H.
\label{3.20}
\end{equation}%
Let $R>0$ be an arbitrary number. Let $\overline{B\left( R\right) }$ be the
closure in the norm of the space $H$ of the open set $B\left( R\right)
\subset H$ of functions $p\left( x,k\right) $ defined as 
\begin{equation}
B\left( R\right) =\left\{ p\in H:\left\Vert p\right\Vert _{H}<R\right\} .
\label{3.21}
\end{equation}

\textbf{Minimization Problem}. \emph{Minimize the functional }$J_{\lambda
}\left( p\right) $\emph{\ on the set }$\overline{B\left( R\right) }.$

\textbf{Remark 3.1}. The analytical part of this paper below is dedicated to
this minimization problem. Since we deal with complex valued functions, we
consider below $J_{\lambda }\left( p\right) $ as the functional with respect
to the 2-D vector of real valued functions $p\left( x,k\right) =\left( %
\mathop{\rm Re}p\left( x,k\right) ,\mathop{\rm Im}p\left( x,k\right) \right)
=\left( p_{1}\left( x,k\right) ,p_{2}\left( x,k\right) \right) \in H.$ Thus,
even though we the consider complex conjugations below, this is done only
for the convenience of writing. Below $\left[ ,\right] $ is the scalar
product in $H$. Even though we use in (\ref{3.16}) and (\ref{3.17}) the
functions $p_{0}=p_{0}\left( k\right) ,$ $p_{1}=p_{1}\left( k\right) ,$ it
is always clear from the context below what do we actually mean in each
particular case: the first component of $p_{1}\left( x,k\right) $ of the
vector function $p\left( x,k\right) $ or the above functions $p_{0}\left(
k\right) ,p_{1}\left( k\right) .$

\section{The Global Strict Convexity of $J_{\protect\lambda }\left( p\right) 
$}

\label{sec:4}

Theorem 4.1 is the main analytical result of this paper.

\textbf{Theorem 4.1}. \emph{\ Assume that conditions of Theorem 3.1 are
satisfied. Then the functional }$J_{\lambda }\left( p\right) $\emph{\ has
the Frech\'{e}t derivative }$J_{\lambda }^{\prime }\left( p\right) $\emph{\
for all }$p\in H.$ \emph{Also, there exists a sufficiently large number }$%
\lambda _{0}=\lambda _{0}\left( r^{\ast },\underline{k},\overline{k}%
,\left\Vert p_{1}\right\Vert _{C\left[ \underline{k},\overline{k}\right]
},R\right) >1$ \emph{depending only on listed parameters and a generic
constant }$C>0$\emph{, such that for all }$\lambda \geq \lambda _{0}$\emph{\
the functional }$J_{\lambda }\left( p\right) $\emph{\ is strictly convex on }%
$\overline{B\left( R\right) },$\emph{\ i.e. for all }$p,p+h\in \overline{%
B\left( R\right) }$\emph{\ }%
\begin{equation}
J_{\lambda }\left( p+h\right) -J_{\lambda }\left( p\right) -J_{\lambda
}^{\prime }\left( p\right) \left( h\right) \geq C\left\Vert h\right\Vert
_{H}^{2}.  \label{3.22}
\end{equation}

\textbf{Proof.} Everywhere below in this paper $C_{2}=C_{2}\left( r^{\ast },%
\underline{k},\overline{k},\left\Vert p_{1}\left( k\right) \right\Vert _{C%
\left[ \underline{k},\overline{k}\right] },R\right) >0$ \ denotes different
constants depending only on listed parameters. Since conditions of Theorem
3.1 are satisfied, then by (\ref{3.15})%
\begin{equation}
\left\Vert V_{\alpha \left( \delta \right) }\right\Vert _{C^{1}\left[ 0,1%
\right] }\leq \left\Vert V^{\ast }\right\Vert _{C^{1}\left[ 0,1\right]
}+C_{1}\delta \leq C_{2}.  \label{3.220}
\end{equation}%
Let $h=\left( h_{1},h_{2}\right) ,$ where $h_{1}=\mathop{\rm Re}h,h_{2}=%
\mathop{\rm Im}h. $ Then (\ref{3.130}), (\ref{3.19}) and (\ref{3.21}) imply
that 
\begin{equation}
\displaystyle\int\limits_{\underline{k}}^{\overline{k}}\left\Vert h\left(
x,k\right) \right\Vert _{C^{1}\left[ 0,1\right] }^{2}dk\leq C_{2}.
\label{3.23}
\end{equation}%
Using (\ref{3.23}), we obtain%
\begin{equation}
\left\vert \displaystyle\int\limits_{k}^{\overline{k}}h^{\prime }\left(
x,k\right) dk\right\vert ^{2}\leq \left( \overline{k}-\underline{k}\right) %
\displaystyle\int\limits_{\underline{k}}^{\overline{k}}\left\vert h^{\prime
}\left( x,k\right) \right\vert ^{2}dk  \label{3.24}
\end{equation}%
\begin{equation*}
\leq \left( \overline{k}-\underline{k}\right) \displaystyle\int\limits_{%
\underline{k}}^{\overline{k}}\left\Vert h\left( x,k\right) \right\Vert
_{C^{1}\left[ 0,1\right] }^{2}dk\leq C_{2}.
\end{equation*}

We use the formula 
\begin{equation}
\left\vert a\right\vert ^{2}-\left\vert b\right\vert ^{2}=\left( a-b\right) 
\overline{a}+\left( \overline{a}-\overline{b}\right) b,\forall a,b\in 
\mathbb{C},  \label{3.250}
\end{equation}%
where $\overline{z}$ is the complex conjugate of $z\in \mathbb{C}$. Denote 
\begin{equation}
a=L\left( p+h\right) ,b=L\left( p\right) ,  \label{3.26}
\end{equation}%
Consider functions $A\left( x,k\right) ,A_{1}\left( x,k\right) ,A_{2}\left(
x,k\right) $ defined as%
\begin{equation}
A=\left\vert L\left( p+h\right) \right\vert ^{2}-\left\vert L\left( p\right)
\right\vert ^{2},A_{1}\left( x,k\right) =\left( a-b\right) \overline{a}%
,A_{2}\left( x,k\right) =\left( \overline{a}-\overline{b}\right) b.
\label{3.25}
\end{equation}%
First, using (\ref{3.17}) and (\ref{3.25}), we single out in $A$ the part,
which is linear with respect to the vector function $h=\left(
h_{1},h_{2}\right) $. Then%
\begin{equation*}
a=\left( p+h\right) ^{\prime \prime }+m^{\prime }-2ik\left( p^{\prime
}+h^{\prime }+m^{\prime }\right)
\end{equation*}%
\begin{equation}
+2k^{2}\left( p^{\prime }+h^{\prime }+m^{\prime }\right) \left(
-\int\limits_{k}^{\overline{k}}\left( p^{\prime }+h^{\prime }+m^{\prime
}\right) \left( x,\tau \right) d\tau +V_{\alpha \left( \delta \right)
}^{\prime }\left( x\right) \right)  \label{3.260}
\end{equation}%
\begin{equation*}
-2i\left( -\int\limits_{k}^{\overline{k}}\left( p^{\prime }+h^{\prime
}+m^{\prime }\right) \left( x,\tau \right) d\tau +V_{\alpha \left( \delta
\right) }^{\prime }\left( x\right) \right) +2k\left( -\int\limits_{k}^{%
\overline{k}}\left( p^{\prime }+h^{\prime }+m^{\prime }\right) \left( x,\tau
\right) d\tau +V_{\alpha \left( \delta \right) }^{\prime }\left( x\right)
\right) ^{2}.
\end{equation*}%
By (\ref{3.25}) 
\begin{equation*}
A_{1}\left( x,k\right) =\left( a-b\right) \overline{a}=\left\{ h^{\prime
\prime }-\left[ 2ik+2k^{2}\left( \int\limits_{k}^{\overline{k}}\left(
p^{\prime }+m^{\prime }\right) \left( x,\tau \right) d\tau +V_{\alpha \left(
\delta \right) }^{\prime }\left( x\right) \right) \right] h^{\prime
}\right\} \overline{a}
\end{equation*}%
\begin{equation}
+\left[ 2k^{2}\left( p^{\prime }+m^{\prime }\right) +2i-4k\left(
\int\limits_{k}^{\overline{k}}\left( p^{\prime }+m^{\prime }\right) \left(
x,\tau \right) d\tau -V_{\alpha \left( \delta \right) }^{\prime }\left(
x\right) \right) \right] \int\limits_{k}^{\overline{k}}h^{\prime }\left(
x,\tau \right) d\tau \cdot \overline{a}  \label{3.261}
\end{equation}%
\begin{equation*}
+\left[ -2k^{2}h^{\prime }\int\limits_{k}^{\overline{k}}h^{\prime }\left(
x,\tau \right) d\tau +2k\left( \int\limits_{k}^{\overline{k}}h^{\prime
}\left( x,\tau \right) d\tau \right) ^{2}\right] \overline{a}.
\end{equation*}%
Hence,%
\begin{equation*}
A_{1}\left( x,k\right) =h^{\prime \prime }\overline{L\left( p\right) }-\left[
2ik+2k^{2}\left( \int\limits_{k}^{\overline{k}}\left( p^{\prime }+m^{\prime
}\right) \left( x,\tau \right) d\tau -V_{\alpha \left( \delta \right)
}^{\prime }\left( x\right) \right) \right] \overline{L\left( p\right) }%
h^{\prime }
\end{equation*}%
\begin{equation}
+\left[ \left( p^{\prime }+m^{\prime }\right) +2i-4k\left( \int\limits_{k}^{%
\overline{k}}\left( p^{\prime }+m^{\prime }\right) \left( x,\tau \right)
d\tau -V_{\alpha \left( \delta \right) }^{\prime }\left( x\right) \right) %
\right] \overline{L\left( p\right) }\int\limits_{k}^{\overline{k}}h^{\prime
}\left( x,\tau \right) d\tau  \label{3.27}
\end{equation}%
\begin{equation*}
+\left\vert h^{\prime \prime }\left( x,k\right) \right\vert ^{2}+\widetilde{A%
}_{1,p}\left( h\right) \left( x,k\right) ,
\end{equation*}%
where $\widetilde{A}_{1,p}\left( h\right) \left( x,k\right) $ depends
nonlinearly on the vector function $\left( h_{1},h_{2}\right) \left(
x,k\right) $. Also, by (\ref{3.220})-(\ref{3.24}) and the Cauchy-Schwarz
inequality%
\begin{equation}
\left\vert \widetilde{A}_{1,p}\left( h\right) \right\vert \left( x,k\right)
\leq \frac{1}{2}\left\vert h^{\prime \prime }\left( x,k\right) \right\vert
^{2}+C_{2}\left( \left\vert h^{\prime }\left( x,k\right) \right\vert ^{2}+%
\displaystyle\int\limits_{\underline{k}}^{\overline{k}}\left\vert h^{\prime
}\left( x,\tau \right) \right\vert ^{2}d\tau \right) .  \label{3.28}
\end{equation}%
To explain the presence of the multiplier \textquotedblleft 1/2" at $%
\left\vert h^{\prime \prime }\left( x,k\right) \right\vert ^{2}$ in (\ref%
{3.28}), we note that it follows from (\ref{3.260}) that the term $h^{\prime
\prime }\overline{a}$ in (\ref{3.261}) contains the term $\left\vert
h^{\prime \prime }\right\vert ^{2},$ which is included in (\ref{3.27})
already, as well as terms%
\begin{equation}
h^{\prime \prime }\overline{h^{\prime }},\text{ }h^{\prime \prime }%
\displaystyle\int\limits_{\underline{k}}^{\overline{k}}\overline{h^{\prime }}%
\left( x,\tau \right) d\tau ,\text{ }h^{\prime \prime }\left( \displaystyle%
\int\limits_{\underline{k}}^{\overline{k}}\overline{h^{\prime }}\left(
x,\tau \right) d\tau \right) ^{2}.  \label{3.280}
\end{equation}%
We now show how do we estimate the third term in (\ref{3.280}), since
estimates of two other terms are simpler. We use the so-called
\textquotedblleft Cauchy-Schwarz inequality with $\varepsilon ",$%
\begin{equation*}
-\left\vert \left( c,d\right) \right\vert \geq -\frac{\varepsilon }{2}%
\left\vert c\right\vert ^{2}-\frac{1}{2\varepsilon }\left\vert d\right\vert
^{2},\forall c,d\in \mathbb{R}^{n},\forall \varepsilon >0,
\end{equation*}%
where $\left( ,\right) $ is the scalar product in $\mathbb{R}^{n}.$ Hence,%
\begin{equation*}
-\left\vert h^{\prime \prime }\left( \displaystyle\int\limits_{k}^{\overline{%
k}}\overline{h^{\prime }}\left( x,\tau \right) d\tau \right) ^{2}\right\vert
\geq -\frac{\varepsilon }{2}\left\vert h^{\prime \prime }\right\vert ^{2}-%
\frac{\left( \overline{k}-\underline{k}\right) }{2\varepsilon }\displaystyle%
\int\limits_{\underline{k}}^{\overline{k}}\left\vert h^{\prime }\left(
x,\tau \right) \right\vert ^{2}d\tau .
\end{equation*}%
Thus, choosing appropriate numbers $\varepsilon >0,$ we obtain the term $%
\left\vert h^{\prime \prime }\left( x,k\right) \right\vert ^{2}/2$ in (\ref%
{3.28}). The second term in the right hand side of (\ref{3.28}) is obtained
similarly.\qquad \qquad\ 

Analogously, using (\ref{3.250})-(\ref{3.25}), we obtain%
\begin{equation*}
A_{2}\left( x,k\right) =\left( \overline{a}-\overline{b}\right) b=\overline{%
\left\{ h^{\prime \prime }-\left[ 2ik+2k^{2}\left( \int\limits_{k}^{%
\overline{k}}\left( p^{\prime }+m^{\prime }\right) \left( x,\tau \right)
d\tau -V_{\alpha \left( \delta \right) }^{\prime }\left( x\right) \right) %
\right] h^{\prime }\right\} }\cdot L\left( p\right)
\end{equation*}%
\begin{equation}
+\overline{\left[ \left( p^{\prime }+m^{\prime }\right) +2i-4k\left(
\int\limits_{k}^{\overline{k}}\left( p^{\prime }+m^{\prime }\right) \left(
x,\tau \right) d\tau -V_{\alpha \left( \delta \right) }^{\prime }\left(
x\right) \right) \right] \int\limits_{k}^{\overline{k}}h^{\prime }\left(
x,\tau \right) d\tau }\cdot L\left( p\right)  \label{3.29}
\end{equation}%
\begin{equation*}
+\widetilde{A}_{2,p}\left( h\right) \left( x,k\right) ,
\end{equation*}%
where $\widetilde{A}_{2,p}\left( h\right) \left( x,k\right) $ depends
nonlinearly on the vector function $h=\left( h_{1},h_{2}\right) $ and
similarly with (\ref{3.28})%
\begin{equation}
\left\vert \widetilde{A}_{2,p}\left( h\right) \right\vert \left( x,k\right)
\leq C_{2}\left( \left\vert h^{\prime }\left( x,k\right) \right\vert ^{2}+%
\displaystyle\int\limits_{\underline{k}}^{\overline{k}}\left\vert h^{\prime
}\left( x,k\right) \right\vert ^{2}dk\right) .  \label{3.30}
\end{equation}

It is clear from (\ref{3.25}), (\ref{3.27})-(\ref{3.30}) that the linear
with respect to the vector function $h=\left( h_{1},h_{2}\right) $ part of $%
A $ consists of the sum of the first two lines of (\ref{3.27}) with the
first two lines of (\ref{3.29}). We denote this linear part as $D_{p}\left(
h\right) \left( x,k\right) .$ Then 
\begin{equation*}
A\left( x,k\right) =A_{1}\left( x,k\right) +A_{2}\left( x,k\right)
=D_{p}\left( h\right) \left( x,k\right) +\widetilde{A}_{1,p}\left( h\right)
\left( x,k\right) +\widetilde{A}_{2,p}\left( h\right) \left( x,k\right) .
\end{equation*}%
Thus, using (\ref{3.20}) and (\ref{3.25}), we obtain%
\begin{equation}
J_{\lambda }\left( p+h\right) -J_{\lambda }\left( p\right) =e^{2\lambda }%
\displaystyle\int\limits_{\underline{k}}^{\overline{k}}\displaystyle%
\int\limits_{0}^{1}D_{p}\left( h\right) \left( x,k\right) e^{-2\lambda x}dxdk
\label{3.31}
\end{equation}%
\begin{equation*}
+e^{2\lambda }\displaystyle\int\limits_{\underline{k}}^{\overline{k}}%
\displaystyle\int\limits_{0}^{1}\left( \widetilde{A}_{1,p}\left( h\right)
\left( x,k\right) +\widetilde{A}_{2,p}\left( h\right) \left( x,k\right)
\right) e^{-2\lambda x}dxdk.
\end{equation*}%
Consider the expression $\widetilde{D}_{p,\lambda }\left( h\right) ,$%
\begin{equation}
\widetilde{D}_{p,\lambda }\left( h\right) =e^{2\lambda }\displaystyle%
\int\limits_{\underline{k}}^{\overline{k}}\displaystyle\int%
\limits_{0}^{1}D_{p}\left( h\right) \left( x,k\right) e^{-2\lambda x}dxdk.
\label{3.32}
\end{equation}%
It follows from (\ref{3.17}), (\ref{3.220}), (\ref{3.27}) and (\ref{3.29})
that $\widetilde{D}_{p,\lambda }\left( h\right) :H\rightarrow \mathbb{R}$ is
a bounded linear functional. Hence, by Riesz theorem, there exists unique
element $M_{p,\lambda }\in H$ such that 
\begin{equation}
\widetilde{D}_{p,\lambda }\left( h\right) =\left[ M_{p,\lambda },h\right]
,\forall h\in H.  \label{3.33}
\end{equation}%
It follows from (\ref{3.28}) and (\ref{3.30})-(\ref{3.33}) that%
\begin{equation*}
J_{\lambda }\left( p+h\right) -J_{\lambda }\left( p\right) -\left[
M_{p,\lambda },h\right] =O\left( \left\Vert h\right\Vert _{H}^{2}\right) .
\end{equation*}%
Thus, the Frech\'{e}t derivative $J_{\lambda }^{\prime }\left( p\right) \in
H $ of the functional $J_{\lambda }\left( p\right) $ at the point $p$ exists
and 
\begin{equation}
J_{\lambda }^{\prime }\left( p\right) =M_{p,\lambda }.  \label{3.34}
\end{equation}

Note that 
\begin{equation}
e^{2\lambda }e^{-2\lambda x}\geq 1,\forall x\in \left[ 0,1\right] .
\label{3.35}
\end{equation}%
Hence, using (\ref{3.28}), (\ref{3.30})-(\ref{3.34}) and Lemma 3.1, we obtain%
\begin{equation*}
J_{\lambda }\left( p+h\right) -J_{\lambda }\left( p\right) -J_{\lambda
}^{\prime }\left( p\right) \left( h\right) \geq
\end{equation*}%
\begin{equation*}
\frac{e^{2\lambda }}{2}\displaystyle\int\limits_{\underline{k}}^{\overline{k}%
}\displaystyle\int\limits_{0}^{1}\left\vert h^{\prime \prime }\right\vert
^{2}\left( x,k\right) e^{-2\lambda x}
\end{equation*}%
\begin{equation}
-C_{2}e^{2\lambda }\displaystyle\int\limits_{\underline{k}}^{\overline{k}}%
\displaystyle\int\limits_{0}^{1}\left\vert h^{\prime }\right\vert ^{2}\left(
x,k\right) e^{-2\lambda x}  \label{3.36}
\end{equation}%
\begin{equation*}
\geq Ce^{2\lambda }\left[ \int\limits_{0}^{1}|h^{\prime \prime
}|^{2}e^{-2\lambda x}dx+\lambda \int\limits_{0}^{1}|h^{\prime
}|^{2}e^{-2\lambda x}dx+\lambda ^{3}\int\limits_{0}^{1}|h|^{2}e^{-2\lambda
x}dx\right]
\end{equation*}%
\begin{equation*}
-C_{2}e^{2\lambda }\displaystyle\int\limits_{\underline{k}}^{\overline{k}}%
\displaystyle\int\limits_{0}^{1}\left\vert h^{\prime }\right\vert ^{2}\left(
x,k\right) e^{-2\lambda x}.
\end{equation*}%
Choose the number $\lambda _{0}=\lambda _{0}\left( r^{\ast },\underline{k},%
\overline{k},\left\Vert p_{1}\left( k\right) \right\Vert _{C\left[ 
\underline{k},\overline{k}\right] },R\right) >1$ so large that $C\lambda
_{0}>2C_{2}.$ Then, using (\ref{3.35}) and (\ref{3.36}), we obtain with a
new generic constant $C>0$ for all $\lambda \geq \lambda _{0}$%
\begin{equation*}
J_{\lambda }\left( p+h\right) -J_{\lambda }\left( p\right) -J_{\lambda
}^{\prime }\left( p\right) \left( h\right) \geq C\left\Vert h\right\Vert
_{H}^{2}.\text{ \ \ \ }\square
\end{equation*}

\section{Global Convergence of the Gradient Projection Method}

\label{sec:5}

Using Theorem 4.1, we establish in this section the global convergence of
the gradient projection method of the minimization of the functional $%
J_{\lambda }\left( p\right) .$ As to some other versions of the gradient
method, they will be discussed in follow up publications.

\subsection{Lipschitz continuity of $J_{\protect\lambda }^{\prime }\left(
p\right) $ with respect to $p$}

\label{sec:5.1}

First, we need to prove the Lipschitz continuity of the functional $%
J_{\lambda }^{\prime }\left( p\right) $ with respect to $p$.

\textbf{Theorem 5.1}.\emph{\ Let conditions of Theorem 3.1 hold. Then the
functional }$J_{\lambda }^{\prime }\left( p\right) $\emph{\ is Lipschitz
continuous on the closed ball }$\overline{B\left( R\right) }.$\emph{\ In
other words,}%
\begin{equation}
\left\Vert J_{\lambda }^{\prime }\left( p^{\left( 1\right) }\right)
-J_{\lambda }^{\prime }\left( p^{\left( 2\right) }\right) \right\Vert
_{H}\leq C_{2}e^{2\lambda }\left\Vert p^{\left( 1\right) }-p^{\left(
2\right) }\right\Vert _{H},\forall p^{\left( 1\right) },p^{\left( 2\right)
}\in \overline{B\left( R\right) },\forall \lambda >0.  \label{5.1}
\end{equation}

\textbf{Proof}. Consider, for example the first line of (\ref{3.27}) for $%
p=p^{\left( 1\right) }$ and denote it $A_{1,1}\left( x,k\right) \left(
p^{\left( 1\right) },h\right) .$ We define $A_{1,1}\left( x,k\right) \left(
p^{\left( 2\right) },h\right) $ similarly. Both these expressions are linear
with respect to $h=\left( h_{1},h_{2}\right) .$ Denote $\widetilde{p}%
=p^{\left( 1\right) }-p^{\left( 2\right) }.$ We have%
\begin{equation*}
A_{1,1}\left( x,k\right) \left( p^{\left( 1\right) },h\right) -A_{1,1}\left(
x,k\right) \left( p^{\left( 2\right) },h\right) =
\end{equation*}%
\begin{equation}
\left( \overline{L\left( p^{\left( 1\right) }\right) }-\overline{L\left(
p^{\left( 2\right) }\right) }\right) \left[ h^{\prime \prime }-\left(
2ik+2k^{2}\left( \int\limits_{k}^{\overline{k}}\left( p^{\left( 1\right)
}+m\right) ^{\prime }\left( x,\tau \right) d\tau -V_{\alpha \left( \delta
\right) }^{\prime }\left( x\right) \right) \right) h^{\prime }\right]
\label{5.2}
\end{equation}%
\begin{equation*}
-2k^{2}\left[ \overline{L\left( p^{\left( 2\right) }\right) }%
\int\limits_{k}^{\overline{k}}\widetilde{p}^{\prime }\left( x,\tau \right)
d\tau \right] h^{\prime }.
\end{equation*}%
It is clear from (\ref{3.17}) that $\left\vert \left( \overline{L\left(
p^{\left( 1\right) }\right) }-\overline{L\left( p^{\left( 2\right) }\right) }%
\right) \right\vert \leq C_{2}\left( \left\vert \widetilde{p}^{\prime \prime
}\right\vert +\left\vert \widetilde{p}^{\prime }\right\vert \right) .$
Hence, using (\ref{3.35}), (\ref{5.2}) and Cauchy-Schwarz inequality, we
obtain%
\begin{equation*}
\left\vert e^{2\lambda }\displaystyle\int\limits_{\underline{k}}^{\overline{k%
}}\displaystyle\int\limits_{0}^{1}\left( A_{1,1}\left( x,k\right) \left(
p^{\left( 1\right) },h\right) -A_{1,1}\left( x,k\right) \left( p^{\left(
2\right) },h\right) \right) e^{-2\lambda x}dxdk\right\vert
\end{equation*}%
\begin{equation*}
\leq C_{2}e^{2\lambda }\left\Vert p^{\left( 1\right) }-p^{\left( 2\right)
}\right\Vert _{H}\left\Vert h\right\Vert _{H}.
\end{equation*}%
The rest of the proof of (\ref{5.1}) is similar. \ $\square $

\subsection{The minimizer of $J_{\protect\lambda }\left( p\right) $ on the
set $\overline{B\left( R\right) }$}

\label{sec:5.2}

Theorem 5.2 claims the existence and uniqueness of the minimizer of the
functional $J_{\lambda }\left( p\right) $ on the set $\overline{B\left(
R\right) }.$

\textbf{Theorem 5.2}. \emph{Let conditions of Theorem 4.1 hold. Then for
every }$\lambda \geq \lambda _{0}$\emph{\ there exists unique minimizer }$%
p_{\min ,\lambda }$\emph{\ of the functional }$J_{\lambda }\left( p\right) $%
\emph{\ on the set }$\overline{B\left( R\right) }.$\emph{\ Furthermore,}%
\begin{equation}
\left[ J^{\prime }\left( p_{\min ,\lambda }\right) ,y-p_{\min ,\lambda }%
\right] \geq 0,\forall y\in \overline{B\left( R\right) }.  \label{5.3}
\end{equation}

\textbf{Proof}. This theorem follows immediately from the above Theorem 4.1
and Lemma 2.1 of \cite{BakKlKosh}. $\square $

Let $Q_{\overline{B}}:H\rightarrow \overline{B\left( R\right) }$ be the
operator of the projection of the space $H$ on the closed ball $\overline{%
B\left( R\right) }.$ Let $\gamma =const.>0$ and let $p^{\left( 0\right) }$
be an arbitrary point of $\overline{B\left( R\right) }$. Consider the
sequence of the gradient projection method,%
\begin{equation}
p^{\left( n+1\right) }=Q_{\overline{B}}\left( p^{\left( n\right) }-\gamma
J_{\lambda }^{\prime }\left( p^{\left( n\right) }\right) \right) ,n=0,1,...
\label{5.4}
\end{equation}

\textbf{Theorem 5.3.} \emph{Let conditions of Theorem 4.1 hold. Then for
every }$\lambda \geq \lambda _{0}$\emph{\ there exists a sufficiently small
number }$\gamma _{0}=\gamma _{0}\left( r^{\ast },\underline{k},\overline{k}%
,\left\Vert p_{0}\right\Vert _{C\left[ \underline{k},\overline{k}\right]
},\left\Vert p_{1}\right\Vert _{C\left[ \underline{k},\overline{k}\right]
},R,\lambda \right) \in \left( 0,1\right) $ \emph{and a number }$q=q\left(
\gamma \right) \in \left( 0,1\right) $\emph{\ such that for every }$\gamma
\in \left( 0,\gamma _{0}\right) $\emph{\ the sequence (\ref{5.4}) converges
to the unique minimizer }$p_{\min ,\lambda }$\emph{\ of the functional }$%
J_{\lambda }\left( p\right) $\emph{\ on the set }$\overline{B\left( R\right) 
}$\emph{\ and }%
\begin{equation}
\left\Vert p^{\left( n\right) }-p_{\min ,\lambda }\right\Vert _{H}\leq
q^{n}\left( \gamma \right) \left\Vert p^{\left( 0\right) }-p_{\min ,\lambda
}\right\Vert _{H},n=1,...  \label{5.5}
\end{equation}

\textbf{Proof}. This theorem follows immediately from the above Theorem 4.1
and Theorem 2.1 of \cite{BakKlKosh}. $\square $

\subsection{Global convergence of the gradient projection method}

\label{sec:5.3}

As it was pointed out in section 3.2, following one of the main concepts of
the regularization theory \cite{BK,Engl}, we assume the existence of the
exact solution $c^{\ast }\left( x\right) $ of our IMSP with the exact, i.e.
noiseless, data $g_{0}^{\ast }\left( k\right) $ in (\ref{2.8}). Below the
superscript \textquotedblleft $^{\ast }$" denotes quantities generated by $%
c^{\ast }\left( x\right) .$ The level of the error $\delta >0$ was
introduced in our data in (\ref{3.14}). In particular, it follows from (\ref%
{3.7}), (\ref{3.70}) and (\ref{3.14}) that%
\begin{equation}
\left\Vert p_{0}-p_{0}^{\ast }\right\Vert _{C\left[ \underline{k},\overline{k%
}\right] },\left\Vert p_{1}-p_{1}^{\ast }\right\Vert _{C\left[ \underline{k},%
\overline{k}\right] }\leq C_{3}\delta ,  \label{5.6}
\end{equation}%
where the number $C_{3}=C_{3}\left( \underline{k},\overline{k}\right) >0$
depends only on listed parameters. Thus, in this section we show that the
gradient projection method delivers points in a small neighborhood of the
function $p^{\ast }$ and, therefore, of the function $c^{\ast }.$ The size
of this neighborhood is proportional to $\delta .$ It is convenient to
indicate in this section dependencies of the functional $J_{\lambda }$ from $%
p_{0},p_{1}$ and $V.$ Hence we write in this section $J_{\lambda }\left(
p,p_{0},p_{1},V_{\alpha \left( \delta \right) }\right) .$

\textbf{Theorem 5.4}. \emph{Assume that conditions of Theorem 4.1 hold.
Also, let the exact function }$p^{\ast }\in B\left( R\right) .$\emph{\ Then
the following accuracy estimates hold for each }$\lambda \geq \lambda _{0}$%
\begin{equation}
\left\Vert p_{\min ,\lambda }-p^{\ast }\right\Vert _{H}\leq C_{2}\delta ,
\label{5.7}
\end{equation}%
\begin{equation}
\left\Vert c_{\min ,\lambda }-c^{\ast }\right\Vert _{L_{2}\left( 0,1\right)
}\leq C_{2}\delta ,  \label{5.8}
\end{equation}%
\emph{where }$p_{\min ,\lambda }$\emph{\ is the minimizer of the functional }%
$J_{\lambda }\left( p,p_{0},p_{1},V_{\alpha \left( \delta \right) }\right) $%
\emph{, which is guaranteed by Theorem 5.2 and }$c_{\min ,\lambda }$\emph{\
is the corresponding reconstructed coefficient (section 3.1). In addition,
let }$\left\{ p^{\left( n\right) }\right\} _{n=0}^{\infty }\subset \overline{%
B\left( R\right) }$\emph{\ be the sequence (\ref{5.4}) of the gradient
projection method, where }$p_{0}$\emph{\ is an arbitrary point of }$%
\overline{B\left( R\right) }$\emph{\ and numbers }$\gamma _{0}$\emph{, }$%
\gamma \in \left( 0,\gamma _{0}\right) $\emph{\ and }$q\left( \gamma \right) 
$\emph{\ are the same as in Theorem 5.3. Let }$\left\{ c_{n}\right\}
_{n=0}^{\infty }$\emph{\ be the corresponding sequence of reconstructed
coefficients (section 3.1). Then the following estimates \ hold}%
\begin{equation}
\left\Vert p^{\left( n\right) }-p^{\ast }\right\Vert _{H}\leq C_{2}\delta
+q^{n}\left( \gamma \right) \left\Vert p_{0}-p_{\min ,\lambda }\right\Vert
_{H},n=1,...,  \label{5.9}
\end{equation}%
\begin{equation}
\left\Vert c_{n}-c^{\ast }\right\Vert _{L_{2}\left( 0,1\right) }\leq
C_{2}\delta +C_{2}q^{n}\left( \gamma \right) \left\Vert p_{0}-p_{\min
,\lambda }\right\Vert _{H},n=1,...  \label{5.10}
\end{equation}

\textbf{Proof}. Obviously%
\begin{equation}
J_{\lambda }\left( p^{\ast },p_{0}^{\ast },p_{1}^{\ast },V^{\ast }\right) =0.
\label{5.11}
\end{equation}%
Using (\ref{3.15}), (\ref{3.170}), (\ref{3.17}), (\ref{5.6}) and (\ref{5.11}%
), we obtain%
\begin{equation*}
J_{\lambda }\left( p^{\ast },p_{0},p_{1},V_{\alpha \left( \delta \right)
}\right) =J_{\lambda }\left( p^{\ast },p_{0}^{\ast }+\left(
p_{0}-p_{0}^{\ast }\right) ,p_{1}^{\ast }+\left( p_{1}-p_{1}^{\ast }\right)
,V^{\ast }+\left( V_{\alpha \left( \delta \right) }-V^{\ast }\right) \right)
\end{equation*}%
\begin{equation}
\leq J_{\lambda }\left( p^{\ast },p_{0}^{\ast },p_{1}^{\ast },V^{\ast
}\right) +C_{2}\left[ \left\vert p_{0}-p_{0}^{\ast }\right\vert +\left\vert
p_{1}-p_{1}^{\ast }\right\vert +\left\vert \left( V_{\alpha \left( \delta
\right) }-V^{\ast }\right) \right\vert \right] J_{\lambda }\left( p^{\ast
}p_{0}^{\ast },p_{1}^{\ast },V^{\ast }\right)  \label{5.12}
\end{equation}%
\begin{equation*}
+C_{2}\left( \left\Vert p_{0}-p_{0}^{\ast }\right\Vert _{C\left[ \underline{k%
},\overline{k}\right] }^{2}+\left\Vert p_{1}-p_{1}^{\ast }\right\Vert _{C%
\left[ \underline{k},\overline{k}\right] }^{2}+\left\Vert V_{\alpha \left(
\delta \right) }-V^{\ast }\right\Vert _{C^{1}\left[ 0,1\right] }^{2}\right)
\leq C_{2}\delta ^{2}.
\end{equation*}%
By Theorems 4.1 and 5.2%
\begin{equation}
J_{\lambda }\left( p^{\ast },p_{0},p_{1},V_{\alpha \left( \delta \right)
}\right) -J_{\lambda }\left( p_{\min ,\lambda },p_{0},p_{1},V_{\alpha \left(
\delta \right) }\right) -\left[ J_{\lambda }^{\prime }\left( p_{\min
,\lambda },p_{0},p_{1},V_{\alpha \left( \delta \right) }\right) ,p^{\ast
}-p_{\min ,\lambda }\right]  \label{5.13}
\end{equation}%
\begin{equation*}
\geq C\left\Vert p_{\min ,\lambda }-p^{\ast }\right\Vert _{H}^{2}.
\end{equation*}%
By (\ref{5.3}) and (\ref{5.12})%
\begin{equation*}
-\left[ J_{\lambda }^{\prime }\left( p_{\min ,\lambda
},p_{0},p_{1},V_{\alpha \left( \delta \right) }\right) ,p^{\ast }-p_{\min
,\lambda }\right] \leq 0,J_{\lambda }\left( p^{\ast },p_{0},p_{1},V_{\alpha
\left( \delta \right) }\right) \leq C_{2}\delta ^{2}.
\end{equation*}
Hence, (\ref{5.13}) implies (\ref{5.7}). Since the function $c_{\min
,\lambda }$ is obtained from the functions $p_{\min ,\lambda }$ and $%
V_{\alpha \left( \delta \right) }$ as described in the end of section 3.1,
then (\ref{5.7}) implies (\ref{5.8}). Next, (\ref{5.9}) follows from (\ref%
{5.5}) and (\ref{5.7}). Finally, (\ref{5.10}) follows from that procedure of
section 3.1 and (\ref{5.8}). $\square $

\textbf{Remark 5.1}. Therefore, Theorem 5.4 ensures the global convergence
property of our method, see the definition in Introduction.

\section{Numerical Studies}

\label{sec:6}

Since the theory of sections 3-5 is the main focus of this paper, we omit
some details of the numerical implementation, both in this and next sections.

\subsection{Algorithm}

\label{sec:6.1}

We now briefly describe our numerical steps for both computationally
simulated and experimental data.

To minimize the functional $J_{\lambda }\left( p\right) ,$ we have written
the derivatives of the operator $L\left( p\right) $ via finite differences
with the step size $h_{x}=0.02$. Also, we have written integrals with
respect to $k$ in discrete forms, using the trapezoidal rule, with the step
size $h_{k}=0.1.$ The differentiation of the data $g_{0}\left( k\right) $
with respect to $k$, which we need in our method (see (\ref{3.70})), was
performed using finite differences with the step size $h_{k}=0.1.$ We have
not observed any instabilities after the differentiation, probably because
the number $h_{k}$ is not too small. Similar conclusions were drawn in works 
\cite{BKSISC, BK, BK2, Chow, KK, Klibfreq, KlibLoc,KTSIAP, Kuzh, IEEE, Exp2, TBKF2}
where similar differentiations were performed, including cases with
experimental data

Next, we have minimized the corresponding discrete version of $J_{\lambda
}\left( p\right) $ with respect to the values of the function $p\left(
x,k\right) $ at those grid points. Initially we have used the gradient
projection method. However, we have observed in our computations that the
regular and simpler gradient method provides practically the same results.
Hence, all computational results below are obtained via the gradient method.
The starting point of this method was $p^{\left( 0\right) }\equiv 0$ and a
specific ball $B\left( R\right) $ was not used. The latter means that
computational results are less pessimistic ones than our theory is. The step
size of the gradient method $\gamma =10^{-5}$ was used. We have observed
that this step size is the optimal one for our computations. The
computations were stopped after 5000 iterations.

Based on our above theory, we have developed the following algorithm:

\begin{enumerate}
\item Find the tail function $V(x)$ via minimizing the functional (\ref{3.12}%
).

\item Minimize the functional (\ref{3.20}). Let $p_{\min ,\lambda }\left(
x,k\right) $ be its minimizer.

\item Calculate the function $q_{\min ,\lambda }\left( x,k\right) =p_{\min
,\lambda }\left( x,k\right) +m\left( x,k\right) ,$ see (\ref{3.16}) and (\ref%
{3.170}).

\item Compute 
\begin{equation*}
v\left( x,\underline{k}\right) =-\int_{\underline{k}}^{\overline{k}}q\left(
x,\tau \right) d\tau +V\left( x\right) .
\end{equation*}

\item Compute the function $\widetilde{c}_{comp}\left( x\right) ,$ see (\ref%
{2.1}) and (\ref{3.5}), 
\begin{equation}
\widetilde{c}_{comp}(x)=\left\vert -v^{\prime \prime }\left( x,\underline{k}%
\right) -\underline{k}^{2}\left( v^{\prime }\left( x,\underline{k}\right)
\right) ^{2}+2i\underline{k}v^{\prime }\left( x,\underline{k}\right)
\right\vert +1.  \label{c}
\end{equation}
\end{enumerate}

In this algorithm, unlike the previous globally convergent algorithms, \cite{BKSISC, BK, Chow, Klibfreq, KlibLoc, Kuzh, IEEE, Exp2, TBKF2}, we do not need to
update the tail function $V(x)$.

%First we solve the problem (\ref{3.10}), (\ref{3.11}) via the minimization of the functional (\ref{3.12}) with condition (\ref{3.13}). Next, having the function $V,$ we minimize the functional $J_{\lambda }\left( p\right) $ in (\ref{3.20}). To minimize functionals (\ref{3.10}) and (\ref{3.20}), we have written differential operators in them via finite differences. Also, integrals with respect to $k$ are naturally written in the discrete forms. Next, we have minimized the discrete form of $J_{\lambda }\left( p\right) $ with respect to the values of unknown discrete function $p\left(x_{j},k_{s}\right) $ at grid points $\left( x_{j},k_{s}\right) $ in both $x$ and $k$ directions. 

\subsection{Numerical testing of computationally simulated data}

\label{sec:6.2}

%To verify the accuracy of our numerical results for experimental data, we have tested our method first on computationally simulated data. 

First, we reconstruct the spatially distributed dielectric constant from
computationally simulated data, which is generated by solving the problem (%
\ref{2.4}), (\ref{2.6}) via the 1-D analog of the Lippmann-Schwinger
equation \cite{KlibLoc}: 
\begin{equation*}
u(x,k)=\frac{\exp (-ik|x-x_{0}|)}{2ik}+k^{2}\int_{0}^{1}\frac{\exp
(-ik|x-\xi |)}{2ik}(c(\xi )-1)u(\xi ,k)\,d\xi .
\end{equation*}%
Here and thereafter, we have use $x_{0}=-1$ in all our computations. Keeping
in mind our desired application to imaging of flash explosive-like targets,
we have chosen in our numerical experiments the true test coefficient $%
c_{true}(x)$ as: 
\begin{equation}
c_{true}(x)=\left\{ 
\begin{array}{cc}
7, & \quad x\in (x_{loc}-d/2,x_{loc}+d/2), \\ 
1, & \quad \mbox{elsewhere},%
\end{array}%
\right.  \label{6.0}
\end{equation}%
where $x_{loc}$ is the location of the center of our target of interest and $%
d$ is its width. Hence, the inclusion/background contrast in (\ref{6.0}) is
7. For our numerical experiments we have chosen in (\ref{6.0}) 
\begin{equation}
x_{loc}=0.1,\,0.2,0.3,0.4\text{ and }d=0.1.  \label{6.2}
\end{equation}

\begin{figure}[tbp]
\begin{center}
\includegraphics[width=0.4\textwidth]{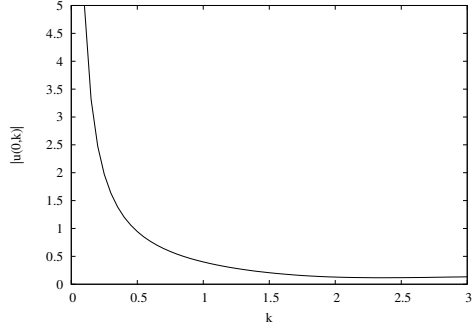}
\end{center}
\caption{The modulus of simulated data on the measurement point $|u(0,k)|$}
\label{fig:u0_abs}
\end{figure}

Figure \ref{fig:u0_abs} displays a typical behavior of the modulus of the
simulated data $\left\vert u\left( 0,k\right) \right\vert $ at the
measurement point $x=0$. One can observe that 
\begin{equation}
\left\vert u\left( 0,k\right) \right\vert \approx 0\text{ for }k>2.
\label{6.20}
\end{equation}
Next, $\left\vert u\left( 0,k\right) \right\vert $ changes too rapidly for $%
k<0.5.$ Hence, the interval $k\in \lbrack 0.5,\,1.5]$ seems to be the
optimal one, and we indeed observed this in our computations. Hence, we
choose for our study $\underline{k}=0.5$ and $\overline{k}=1.5$. We note
that even though the above theory of the choice of the tail function $V(x)$
works only for sufficiently large values of $\overline{k},$ the notion
\textquotedblleft sufficiently large" is relative, see, e.g. (\ref{6.20}).
Besides, it is clear from section 7 that we actually work in the Gigahertz
range of frequencies, and this can be considered as the range of large
frequencies in Physics.

%Figure 1: \emph{A typical behavior of the modulus }$\left\vert w\left( 0,k\right) \right\vert $\emph{\ of the simulated data }$w\left( 0,k\right) .$ \emph{\ One can observe that }$\left\vert w\left( 0,k\right) \right\vert\approx 0$\emph{\ for }$k>???$\emph{\ and it has maximal values near }$k=???$

Next, having the values of $u(0,k)$, we calculate the function $g_{0}(k)$ in
(\ref{2.8}) and introduce the random noise in this function 
\begin{equation*}
g_{0,\,noise}(k)=g_{0}(k)(1.0+0.05\,\sigma (k)),\quad \sigma (k)=\sigma
_{1}(k)+i\sigma _{2}(k),
\end{equation*}%
where $\sigma _{1}(k)$ and $\sigma _{2}(k)$ are random numbers, uniformly
distributed on $(-1,1)$.

The next important question is about the choice of an optimal parameter $%
\lambda =\lambda _{opt}.$ Indeed, even though Theorem 4.1 says that the
functional $J_{\lambda }\left( p\right) $ is strictly convex on the closed
ball $\overline{B\left( R\right) }$ for all $\lambda \geq \lambda _{0},$ in
fact, the larger $\lambda $ is, the less is the influence on $J_{\lambda
}\left( p\right) $ of those points $x\in \left( 0,1\right) ,$ which are
relatively far from the point $\left\{ x=0\right\} $ \ where the data are
given. Hence, we need to choose such a value of $\lambda _{opt},$ which
would provide us satisfactory images of inclusions, whose centers $x_{loc}$
are as in (\ref{6.2}): $x_{loc}\in \left[ 0.1,0.4\right] $.

\begin{figure}[tbp]
\begin{center}
\includegraphics[width=0.4\textwidth]{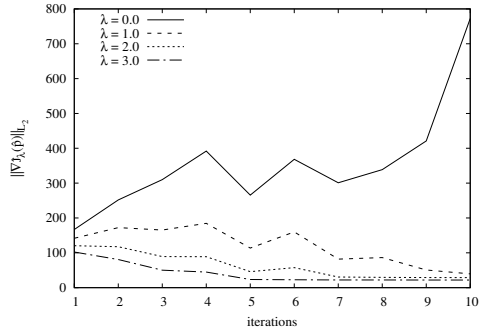}
\end{center}
\caption{The discrete $L_{2}$ norm of the gradient of the functional $J_{%
\protect\lambda }\left( p\right) $ for different $\protect\lambda$}
\label{fig:gnorm}
\end{figure}

Let $\left\Vert \nabla J_{\lambda }\left( p\right) \right\Vert _{L_{2}\left(
\left( 0,1\right) \times \left( \underline{k},\overline{k}\right) \right) }$
be the discrete $L_{2}\left( \left( 0,1\right) \times \left( \underline{k},%
\overline{k}\right) \right) -$ norm of the gradient of the above described
discrete version of the functional $J_{\lambda }\left( p\right) .$ Figure %
\ref{fig:gnorm} displays the dependencies of this norm on the number of
iteration of the gradient method for different values of $\lambda $. We have
observed in our computations that these dependencies are very similar for
targets satisfying (\ref{6.0}), (\ref{6.2}) with different values of
target/background contrasts. One can see that the process diverges at $%
\lambda =0$, which is to be expected, since convexity of $J_{\lambda
=0}\left( p\right) $ is not guaranteed. Also, we observe that the larger $%
\lambda $ is, the faster the process converges. We have found that the
optimal value of $\lambda $ for targets satisfying (\ref{6.2}) is $\lambda
_{opt}=3$.

We also apply a post-processing procedure after step 5 of the above
algorithm. More precisely, we smooth out the function $\widetilde{c}%
_{comp}(x)$ (\ref{c}) using a simple averaging procedure over two
neighboring grid points. Next, the resulting function $\widehat{c}%
_{comp}\left( x\right) $ is truncated as 
\begin{equation}
c_{comp}\left( x\right) =\left\{ 
\begin{array}{ll}
\widehat{c}_{comp}\left( x\right) , & \text{ if }\widehat{c}_{comp}\left(
x\right) \geq 0.8\max (\widehat{c}_{comp}\left( x\right) ), \\ 
1, & \text{ otherwise.}%
\end{array}%
\right.  \label{6.1}
\end{equation}%
The function $c_{comp}\left( x\right) $ in (\ref{6.1}) is considered as our
reconstructed coefficient $c\left( x\right) .$

\begin{figure}[tbp]
\begin{center}
\subfloat[\label{fig:c1}]{\includegraphics[width=0.4
\textwidth]{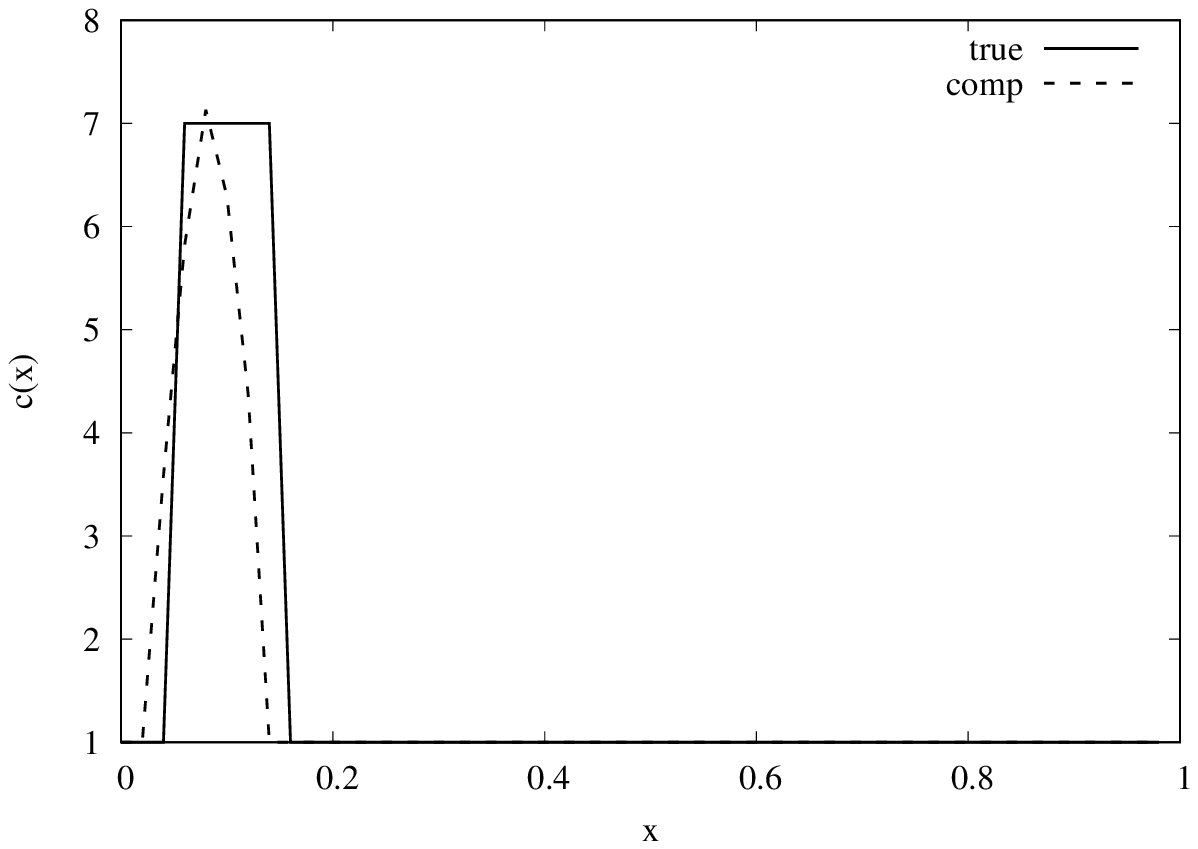}} 
\subfloat[\label{fig:c2}]{
{\includegraphics[width=0.4\textwidth]{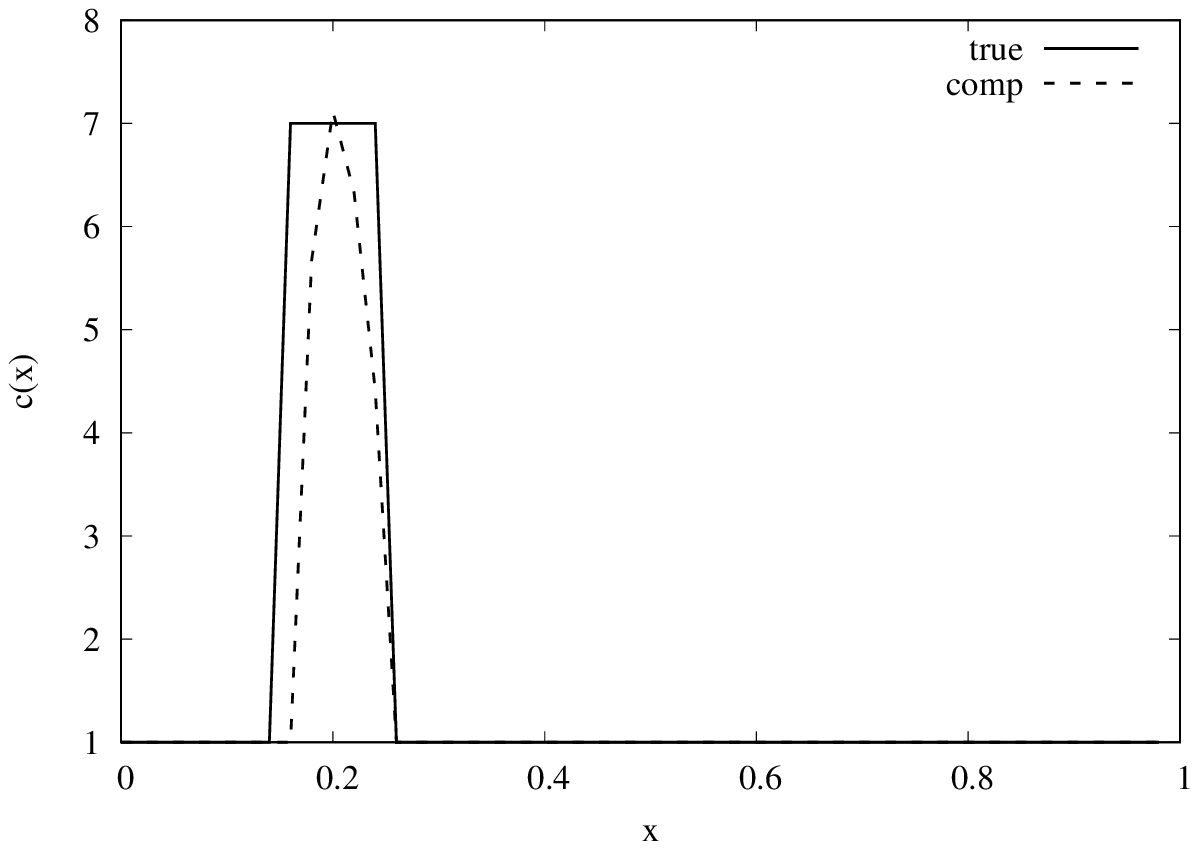}} } \\[0pt]
\subfloat[\label{fig:c3}]{
{\includegraphics[width=0.4\textwidth]{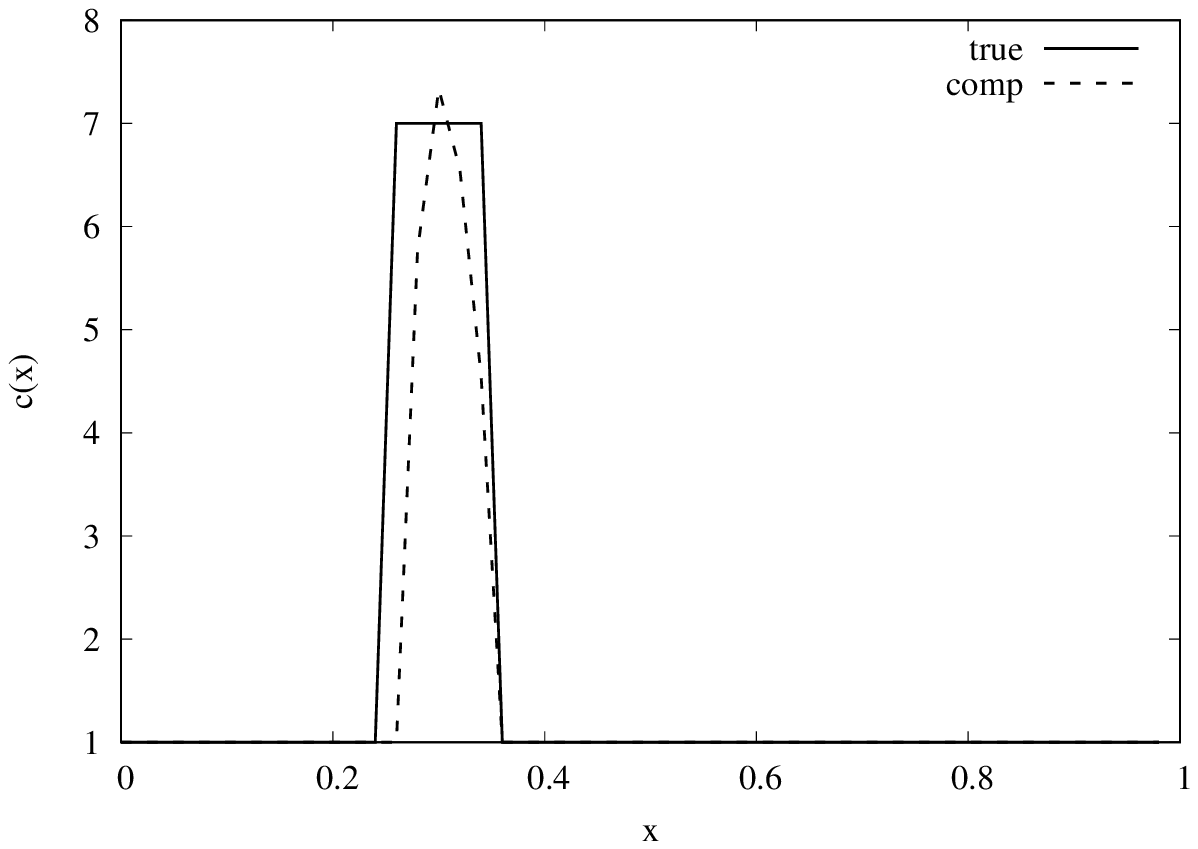}} } 
\subfloat[\label{fig:c4}]{
{\includegraphics[width=0.4\textwidth]{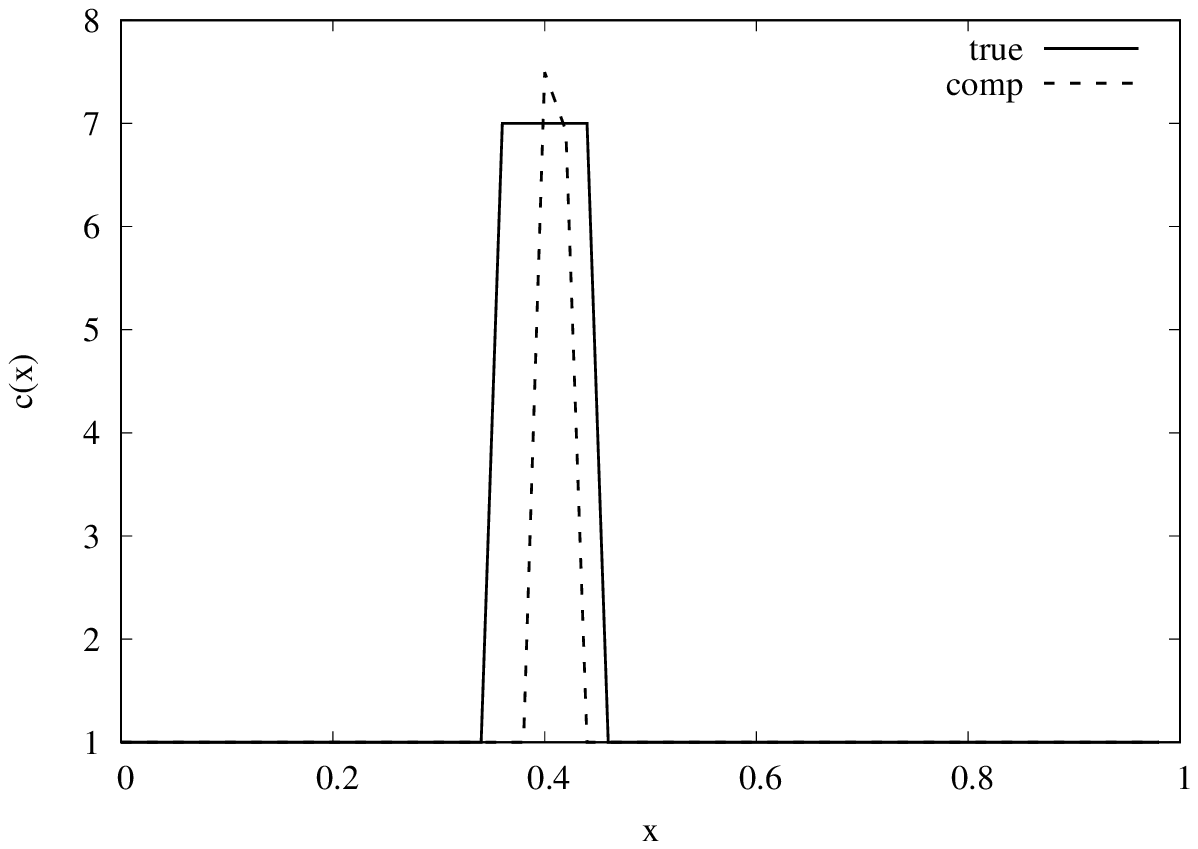}} }
\end{center}
\par
%: (a) picture, (b) scheme}
\caption{True (solid line) and computed (dashed line) coefficients $c(x)$.
a) $x_{loc} = 0.1$, b) $x_{loc} = 0.2$, c) $x_{loc} = 0.3$, d) $x_{loc} =
0.4 $. The optimal value $\protect\lambda _{opt}=3$ was used in all
computations. }
\label{fig:results}
\end{figure}

The computational results $c_{comp}$ for different values of $x_{loc}$ are
shown in Figure \ref{fig:results}. One can see that the proposed algorithm
accurately reconstructs both locations and values of the coefficient $%
c_{true}(x)$. Similar accuracy was obtained for other target/background
contrasts in (\ref{6.0}) varying from 2 to 10.

\section{Numerical Results for Experimental Data}

\label{sec:7}

We use here the same experimental data as ones used in \cite%
{KlibLoc,Kuzh,IEEE}, where these data were treated by the tail functions
method. Thus, it is worth to test the new method of this paper on the same
data set. In \cite{Kuzh,IEEE} the wave propagation process was modeled by a
1-D hyperbolic equation, the Laplace transform with respect to time was
applied to the solution of this equation and then the tail functions method
was applied to the corresponding IMSP. In \cite{KlibLoc} the process was
modeled by IMSP (\ref{2.8}) and the tail functions method was applied to
this IMSP. The data in \cite{Kuzh,IEEE} and in \cite{KlibLoc} were obtained
after applying Laplace and Fourier transforms respectively to the original
time dependent data.

We have observed a substantial mismatch of amplitudes between
computationally simulated and experimental data. Hence, we have calibrated
experimental data here via multiplying them by the calibration factor $%
10^{-7},$ just as in \cite{KlibLoc,Kuzh,IEEE}.

\subsection{Data collection}

\label{sec:7.1}

\begin{figure}[tbp]
\begin{center}
\begin{tikzpicture}[font=\footnotesize,scale=0.3]

%a car
\draw [fill=black!20!white,line width=1.2]  (11,0.3) -- (10.5,2) -- (7,2.5) arc (270:200:0.5) -- (5.5,4.5) arc (20:90:0.8)  -- (-1,5) arc (90:180:0.8) -- (-3,0.3) --  (-2.7,0) -- (-1.2,0) arc (180:0:1.2) -- (6.8,0) arc (180:0:1.2) -- (10.7,0) -- (11,0.3) ;
% wheels
\draw [fill=black!50!white,line width=1.2] (0,0) circle (1);
\draw [fill=black!50!white,line width=1.2] (8,0) circle (1);

\draw [fill=black!30!white,line width=1.2] (4.3,5) -- (4,6) -- (3.7,5) -- (4.3,5);

\draw [fill=black!30!white,line width=1.2] (5,6.5) -- (5,7.5) -- (4,7.3) -- (2, 6.3) -- (2, 5.7) --   (3,5.5) -- (5,6.5);

\draw [line width=1.2] (3,5.5) -- (3,6.5) -- (5,7.5);
\draw [line width=1.2] (3,6.5) -- (2,6.3);

\draw [fill=black!70!white,line width=1.2] (3.5,6.2) circle (0.2);
\draw [fill=black!70!white,line width=1.2] (4,6.5) circle (0.2);
\draw [fill=black!70!white,line width=1.2] (4.5,6.8) circle (0.2);

\draw [fill=black!10!white,line width=1.0] (-5,-1) rectangle (30,-6);

\draw [fill=black!50!white,line width=1.0] (20,-1.5) rectangle (25,-3);

\node  at (22.5,-2.25) {Target}; 

\draw (5.5,6) arc (300:350:0.7);
\draw (7.5,5.0) arc (300:350:1.5);
\draw (9.5,4) arc (300:350:2.5);
\draw (11.5,3) arc (300:350:3.5);
\draw (13.5,2) arc (300:350:4.5);
\draw (15.5,1) arc (300:350:5.5);
\draw (17.5,0) arc (300:350:6.5);
\draw (19.5,-1) arc (300:350:7.5);

\end{tikzpicture}
\end{center}
\caption{Schematic diagram of data collection by the Forward Looking Radar
of the US Army Research Laboratory}
\label{fig:setup}
\end{figure}
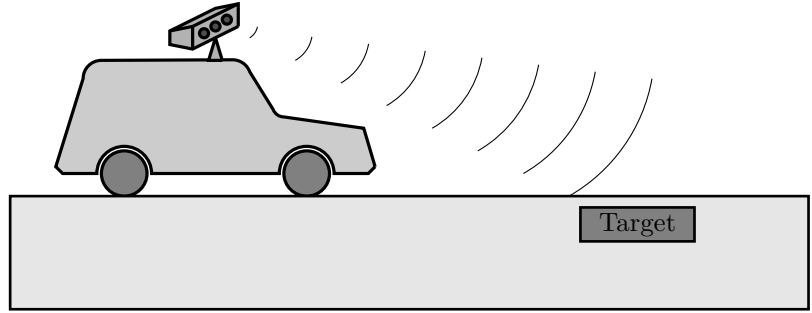

Our experimental data were collected in the field by the Forward Looking
Radar of the US Army Research Laboratory \cite{Radar}. The schematic diagram
of data collection is presented on Figure \ref{fig:setup}. The device has
two sources placed on the top of a car. Sources emit pulses. The device also
has 16 detectors. Detectors measure backscattering time resolved signal,
which is actually the voltage. Pulses of only one component of the electric
field are emitted and the same component is measured on those detectors. The
time step size of measurements is 0.133 nanosecond and the maximal
amplitudes of the measured signal are seen about 2 nanoseconds, see Figure
5. Since 1 nanosecond corresponds to the frequency of 1 Gigahertz \cite{nano}%
, then the corresponding frequency range is in Gigahertz, which are
considered as high frequencies in Physics. The car moves and the time
dependent backscattering signal is measured on distances from 20 to 8 meters
from the target of interest. The collected signals are averaged. Users know
horizontal coordinates of each target with a very good precision: to do
this, the Ground Positioning System is used. Two kinds of targets were
tested: ones located in air and ones buried on the depth of a few
centimeters in the ground.

\begin{figure}[tbp]
\begin{center}
\includegraphics[width=0.4\textwidth]{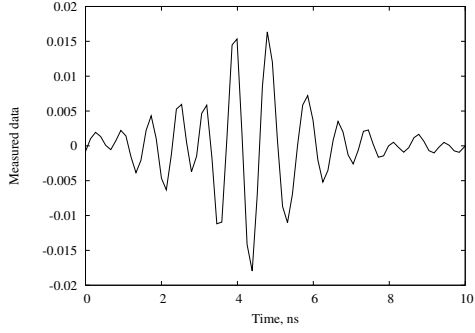}
\end{center}
\caption{Measured time dependent data for bush after being multiplied by the
calibration factor $10^{-7}$. The horizontal axis is time in nanoseconds.}
\label{fig:exp_data}
\end{figure}

\subsection{Results}

\label{sec:7.2}

While it is assumed both in (\ref{2.1}) and (\ref{6.1}) that $c\left(
x\right) \geq 1,$ we had one target buried in the ground, in which $0<c<1.$
This target was a plastic cylinder. It was shown on page 2944 of \cite{IEEE}
that, using the original time dependent date, one can figure out that inside
the target $c\in \left( 0,1\right) .$ Hence, in this case we replace (\ref{c}%
) and (\ref{6.1}) with%
\begin{equation}
\widetilde{c}_{comp}(x)=1-\left\vert -v^{\prime \prime }\left( x,\underline{k%
}\right) -\underline{k}^{2}\left( v^{\prime }\left( x,\underline{k}\right)
\right) ^{2}+2i\underline{k}v^{\prime }\left( x,\underline{k}\right)
\right\vert ,  \label{7.20}
\end{equation}%
\begin{equation}
c_{comp}\left( x\right) =\left\{ 
\begin{array}{ll}
\widetilde{c}_{comp}\left( x\right) , & \text{ if }\widetilde{c}%
_{comp}\left( x\right) \in \left( 0.1,1\right) , \\ 
1, & \text{ otherwise.}%
\end{array}%
\right.  \label{7.2}
\end{equation}

Suppose that a target occupies a subinterval $I\subset \left( 0,1\right) .$
In fact, we estimate here the ratio of dielectric constants of targets and
backgrounds for $x\in I$. Thus, actually our computed function $%
c_{comp}\left( x\right) $ in (\ref{6.1}) and (\ref{7.2}) is an estimate of
the function $P\left( x\right) ,$%
\begin{equation}
P\left( x\right) =\frac{c_{\text{target}}\left( x\right) }{c_{bckgr}}\approx
c_{comp}\left( x\right) ,x\in I,  \label{7.3}
\end{equation}%
where $c_{\text{target}}\left( x\right) $ is the spatially distributed
dielectric constant of that target. Using (\ref{6.1}), (\ref{7.20}), (\ref%
{7.2}) and (\ref{7.3}),\ we define the computed target/background contrast
in the dielectric constant as%
\begin{equation}
\widetilde{P}=\left\{ 
\begin{array}{c}
\max c_{comp}\left( x\right) ,\text{ if }c_{comp}\left( x\right) \geq
1,\forall x\in \left[ 0,1\right] , \\ 
\min c_{comp}\left( x\right) ,\text{ if }c_{comp}\left( x\right) \leq
1,\forall x\in \left[ 0,1\right] .%
\end{array}%
\right.  \label{7.4}
\end{equation}%
Finally, we introduce the number $c_{est},$ which is our estimate of the
dielectric constant of a target, 
\begin{equation}
c_{est}=c_{bckgr}\widetilde{P}.  \label{7.5}
\end{equation}

We have chosen the interval $\left[ \underline{k},\overline{k}\right] $ as 
\begin{equation}
k\in \left[ 2.7,3.2\right] =\left[ \underline{k},\overline{k}\right] .
\label{7.1}
\end{equation}%
The considerations for the choice (\ref{7.1}) were similar with ones for the
case of simulated data in section 6.2.

\begin{figure}[tbp]
\begin{center}
\subfloat[\label{fig:c_bush}]{\includegraphics[width=0.33
\textwidth]{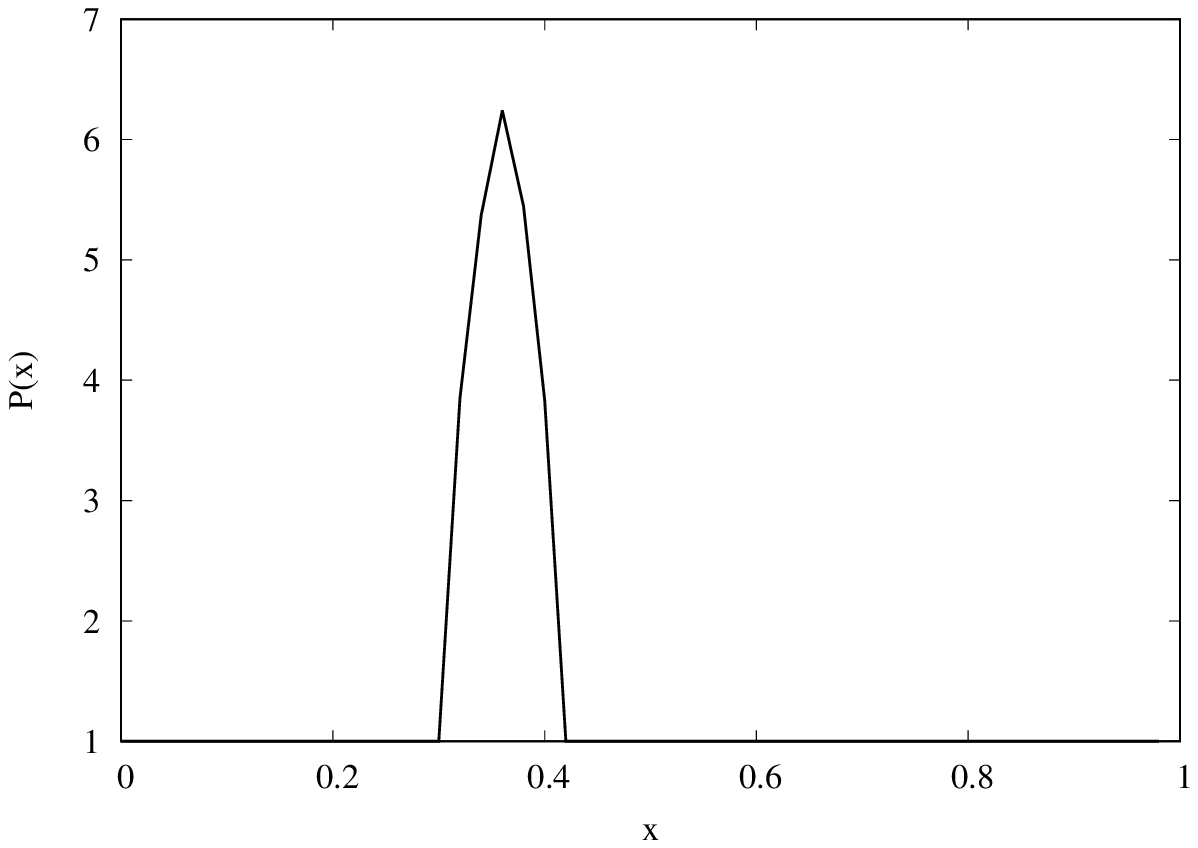}} 
\subfloat[\label{fig:c_mcyl}]{
{\includegraphics[width=0.33\textwidth]{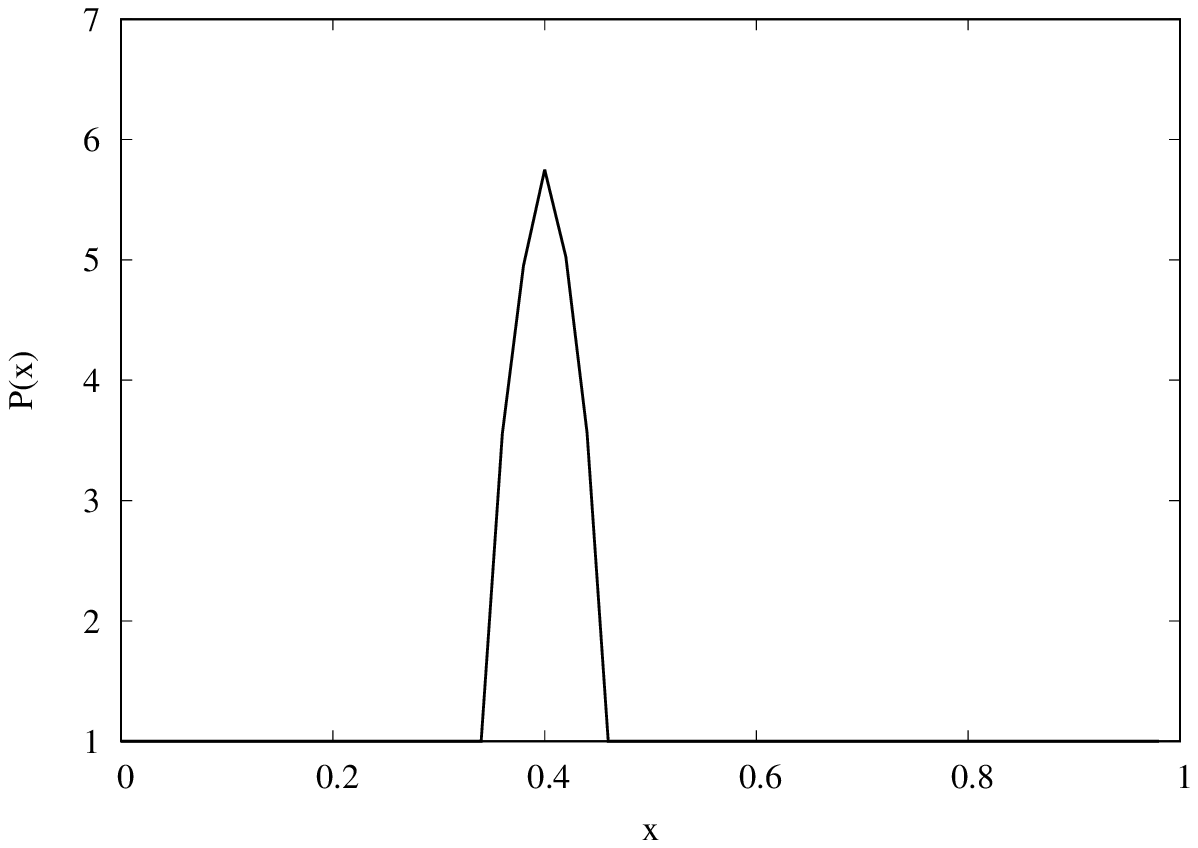}} } 
\subfloat[\label{fig:c_plastic}]{
{\includegraphics[width=0.33\textwidth]{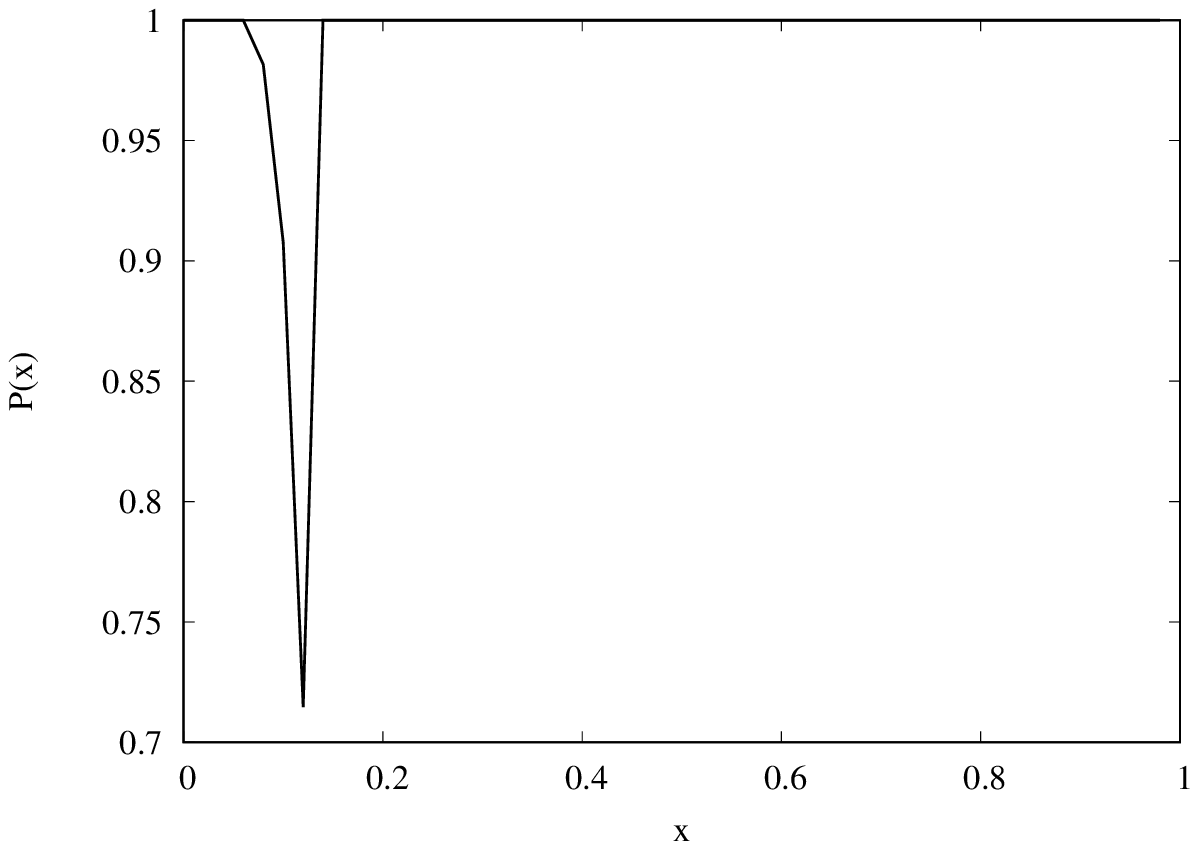}} }
\end{center}
\caption{Computed functions $c_{comp}\left( x\right)$ for different targets:
a) bush in air, b) buried metal cylinder, c) buried plastic cylinder}
\label{fig:exp_res}
\end{figure}

We had experimental data for total five targets. The background was air in
the case of targets placed in air with $c_{bckgr}=1$ and it was sand with $%
c_{bckgr}\in \left[ 3,5\right] $ \cite{Diel} in the case of buried targets.
Two targets, bush and wood stake, were placed in air and three targets,
metal box, metal cylinder and plastic cylinder, were buried in sand. Figures %
\ref{fig:exp_res} display some samples of calculated images of targets.

Dielectric constants of targets were not measured in experiments. So, the
maximum what we can do at this point is to compare our computed values of $%
c_{est}$ with published ones. This is done in Table \ref{tab1}, in which $%
c_{true}$ is a published value. As to the metallic targets, it was
established numerically in \cite{Kuzh,IEEE} that they can be approximated as
dielectric targets with large values of the dielectric constant,%
\begin{equation}
c\in \left[ 10,30\right] .
\end{equation}%
Published values of dielectric constants of sand, wood and plastic can be
found in \cite{Diel}. As to the case when the target was a bush, we took the
interval of published values from \cite{Chuah}. Bush was the most
challenging target to image. This is because bush is obviously a
significantly heterogeneous target.

\begin{table}[!h]
\begin{center}
\begin{tabular}{|l|c|c|c|c|c|}
\hline
Target & buried/no & $\widetilde{P}$ & $c_{bckgr}$ & $c_{est}$ & $c_{true}$
\\ \hline
Bush & no & 6.24 & 1 & 6.24 & $[3,20]$ \\ 
Wood stake & no & 5.43 & 1 & 5.43 & $[2,6]$ \\ 
Metal box & buried & 5.75 & $\left[ 3,5\right] $ & $\left[ 17.25,28.75\right]
$ & $[10,30]$ \\ 
Metal cylinder & buried & 6.48 & $\left[ 3,5\right] $ & $\left[ 19.44,32.40%
\right] $ & $[10,30]$ \\ 
Plastic cylinder & buried & 0.71 & $\left[ 3,5\right] $ & $\left[ 2.13,3.55%
\right] $ & $[1.1,3.2]$ \\ \hline
\end{tabular}%
\end{center}
\caption{Summary of estimated dielectric constants $c_{est}$. }
\label{tab1}
\end{table}

For the engineering part of this team of coauthors (LN and AS), the depth of
burial of a target is not of an interest here since all depths are a few
centimeters.\ It is also clear that it is impossible to figure out the shape
of the target, given so limited information content. On the other hand, the
most valuable piece of the information for LN and AS is in estimates of the
dielectric constants of targets. Therefore, Table \ref{tab1} is the most
interesting piece of the information from the engineering standpoint.
Indeed, one can see in this table that values of estimated dielectric
constants $c_{est}$ are always within limits of $c_{true}.$ As it was
pointed out in section 1, these estimates, even if not perfectly accurate,
can be potentially very useful for the quite important goal of reducing the
false alarm rate. This indicates that the technique of the current paper
might potentially be quite valuable for the goal of an improvement of the
false alarm rate. The above results inspire LN and AS to measure dielectric
constants of targets in the future experiments. Our team plans to treat
those future experimental data by the numerical method of this publication.

\section{Concluding Remarks}

\label{sec:8}

We have developed a new globally convergent numerical method for the 1-D
Inverse Medium Scattering Problem (\ref{2.8}). Unlike the tail function
method, the one of this paper does not impose the smallness condition on the
size of the interval $\left[ \underline{k},\overline{k}\right] $ of wave
numbers. The method is based on the construction of a weighted cost
functional with the Carleman Weight Function in it. The main new theoretical
result of this paper is Theorem 4.1, which claims the strict convexity of
this functional on any closed ball $\overline{B\left( R\right) }\subset H$
for any radius $R>0$, as long as the parameter $\lambda >0$ of this
functional is chosen appropriately. Global convergence of the gradient
method of the minimization of this functional to the exact \ solution is
proved. Numerical testing of this method on both computationally simulated
and experimental data shows good results.

%\bibliographystyle{siamplain}
%\bibliography{biblio}

\end{document}